\documentclass[%
amsmath,amssymb,
preprint,%
]{revtex4-1}

\usepackage{graphicx}
\usepackage{dcolumn}
\usepackage{bm}

\usepackage[utf8]{inputenc}
\usepackage[T1]{fontenc}
\usepackage{mathptmx}
\usepackage{etoolbox}
\usepackage{amssymb,amsmath,amsthm}
\usepackage{mathrsfs}
\usepackage{hyperref}
\usepackage{bm}
\usepackage{color}
\usepackage{makecell}
\usepackage{float}
\usepackage{overpic}
\usepackage{subfig}
\usepackage{threeparttable}

\newtheorem{theorem}{Theorem}[section]
\newtheorem{corollary}[theorem]{Corollary}

\newtheorem{remark}[theorem]{Remark}

\makeatletter
\def\@email#1#2{%
 \endgroup
 \patchcmd{\titleblock@produce}
  {\frontmatter@RRAPformat}
  {\frontmatter@RRAPformat{\produce@RRAP{*#1\href{mailto:#2}{#2}}}\frontmatter@RRAPformat}
  {}{}
}%
\makeatother
\begin{document}

\preprint{ }

\title[Spatiotemporal Patterns on a Circular Domain]{Spatiotemporal Patterns Induced by Turing-Hopf Interaction and Symmetry on a Disk}

\author{Yaqi Chen}
\affiliation{Department of Mathematics, Harbin Institute of Technology, Weihai, Shandong 264209, P.R.China.}

\author{Xianyi Zeng}
\affiliation{Department of Mathematics, Lehigh University, Bethlehem, PA 18015, United States.}

\author{Ben Niu~*}%
 \email{niu@hit.edu.cn.}
\affiliation{Department of Mathematics, Harbin Institute of Technology, Weihai, Shandong 264209, P.R.China. 
}%

\date{\today}

\begin{abstract}
Turing bifurcation and Hopf bifurcation are two important kinds of transitions giving birth to inhomogeneous solutions, in spatial or temporal ways. On a disk, these two bifurcations may lead to equivariant Turing-Hopf bifurcations. In this paper, normal forms for three kinds of Turing-Hopf bifurcations are given and the breathing, standing wave-like, and rotating wave-like patterns are found in numerical examples.
\end{abstract}

\maketitle

\section{\label{sec1}Introduction}

In many reaction-diffusion equations, complex Turing-Hopf spatiotemporal patterns may appear, which is usually used as an explanation for many dynamic phenomena in chemical reactions, epidemics, and metapopulation models \citep{Camara2016J,Yang2016J,Cao2018J,Kumari2020J,Kumar2020J}.
Mathematically, the solution of the reaction-diffusion equation is often written in the form of $u_j(\vec{x},t)=u_{j0}\mathrm{e}^{\mathrm{i}\vec{q}\cdot\vec{x}+\sigma t}$ \citep{Cross1993J} to characterize its spatiotemporal dynamics, where $\vec{q}$ is the wave vector and $\sigma$ is the eigenvalue with the largest real part. For a Turing-Hopf bifurcation,  $\vec{q}$ is nonzero and $\sigma$ is also an imaginary value $\mathrm{i}\omega$. Thus, there exists the interaction of two Fourier modes \citep{Perraud1993J,Heidemann1993J,Vallette1994J,Song2014J,An2018J}.

There are many methods for studying the Turing-Hopf bifurcation.
In numerical experiments, simple hexagonal arrangements of spots, stripes, and other spatiotemporal patterns can be observed in chemical reactions, semiconductors, predator-prey models, and other systems through numerical tools, indicating that this bifurcation can reveal some complex spatiotemporal patterns of dynamic systems \citep{Rovinsky1992J,Meixner1997J,Bose2000J,Baurmann2007J}.
In physics, researchers often use multi-scale methods to derive the amplitude equation of the Turing-Hopf bifurcation to analyze pattern formations \citep{Just2001J,Venkov2007J,Ledesma2020J}.
Moreover, in recent years, scholars have begun to use normal forms to analyze Turing-Hopf bifurcation.
In particular, Song et al. \citep{Song2016J} and Jiang et al. \citep{Jiang2020J} derived the normal form of the Turing-Hopf bifurcation of partial differential equations (PDEs) and partial functional differential equations (PFDEs), respectively.
Following the method proposed, there are many subsequent works on normal forms of the Turing-Hopf bifurcation \citep{Song2019J,Wu2020J,Lv2021J,Duan2022J}.

However, most of these previous works have focused on one-dimensional intervals, and a few of them considered two-dimensional spaces, for example, rectangular domains  \citep{Just2001J,Cao2018J,Chen2021J} and circular domains\citep{Abid2015J,Paquin-Lefebvre2019J}, in which more abundant spatiotemporal patterns will be generated.
In fact, the complex spatiotemporal patterns appearing in circular domains can be studied through the equivariant bifurcation. Precisely, the works on symmetric group theory in \citep{Golubitsky1989M} and the theory of equivariant normal forms in \citep{Guo2013M} are required. Equivariant Turing-Hopf bifurcation with time delay on a disk has not been considered, to our best knowledge. Therefore, in this paper, we consider a delayed general reaction-diffusion system with homogeneous Neumann boundary conditions on a disk and aim to explain many interesting spatiotemporal patterns induced by Turing-Hopf interaction and the symmetry.

Compared to previous work, our research in this paper has several new features. We gave formulas of the equivariant normal forms truncated to the third order of a general reaction-diffusion system on a disk and divided them into three types: \textit{ET-H}, \textit{T-EH}, and \textit{ET-EH bifurcations}, according to different structure of the center subspace of the equilibrium.
We characterize the long-term asymptotic behavior of the solution by normal forms, which can explain the occurrence of many patterns in real life more fitly. The theoretical results indicate the existence of several kinds of interesting patterns, including \textit{breathing}, \textit{standing wave-like}, \textit{rotating wave-like patterns} and so on.

The rest of the paper is organized as follows. In Section \ref{sec-existence}, We give preliminaries required for normal form derivation, including the definition of phase space, the eigenvalue problem of the Laplace operator on a circular domain, and the necessary assumptions for bifurcation. In Section \ref{sec3-Normal}, main results of normal forms for \textit{ET-H}, \textit{T-EH}, and \textit{ET-EH bifurcations} on a disk are shown respectively. In Section \ref{sec4-Numerical Simulations},
two delayed mussel-algae systems are selected. Rich spatiotemporal patterns are observed near the Turing-Hopf points.

\section{Preliminaries}\label{sec-existence}

We consider a general delayed reaction-diffusion system of $n$ equations with homogeneous Neumann boundary conditions defined on a disk as follows:
\begin{equation}\label{general reaction-diffusion equations r theta}
\frac{\partial U(t,r,\theta)}{\partial t}=D(\mu)\Delta_{r\theta} U(t,r,\theta)+L(\mu)U_t(r,\theta)+F\left(U_t(r,\theta),\mu\right),~(r,\theta)\in \mathbb{D},~t>0,
\end{equation}
where $\Delta_{r \theta} U=\frac{\partial^{2}U}{\partial r^{2}} +\frac{1}{r} \cdot \frac{\partial U}{\partial r}+\frac{1}{r^{2}} \cdot \frac{\partial^{2} U}{\partial \theta^{2}},~\mathbb{D}=\{(r, \theta): 0 \leq r \leq R, 0 \leq \theta \leq 2 \pi\},~{U}_{t} (\vartheta)(r,\theta)=U(t+\vartheta,r,\theta),~\vartheta \in [-1,0)$. Here, we normalize the maximum delay to 1. To study the interaction between Turing instability and Hopf bifurcation, we usually select two parameters, i.e. $\mu=(\mu_1,\mu_2) \in \mathbb{R}^2$.

When considering a reaction-diffusion equation with time delay, one usually use the phase space of functions \citep{Hale1977M,Wu1996M},
$$
{\mathscr{C}}:=C([-1,0], {\mathscr{X}_{\mathbb{C}}}),
$$
where ${\mathscr{X}_{\mathbb{C}}}$ is the complexification of
$$
{\mathscr{X}}=\left\{\tilde{u}(r, \theta)\in {W}^{2,2}(\mathbb{D}): \partial_r \tilde{u}(R, \theta)=0,~\theta \in [0,2\pi) \right\},
$$
with $\mathbb{L}^2$ inner product $\langle u(r,\theta),v(r,\theta)\rangle=\iint_{\mathbb{D}}r u(r,\theta) \bar{v}(r,\theta) \mathrm{d} r\mathrm{d}\theta$ weighted $r$.
${U}_{t} \in \mathscr{C}^n, ~L:\mathbb{R}\times {\mathscr{C}^n}\rightarrow\mathscr{X}_{\mathbb{C}}^n$ is a bounded linear operator, and
${F}: {\mathscr{C}}^n \times \mathbb{R} \rightarrow \mathscr{X}_{\mathbb{C}}^n$ is $C^k~(k\ge 3)$. Here we only consider the zero equilibrium, that is to say, we assume $F(0,\mu)=0,~\forall \mu \in \mathbb{R}^2$.

For $\mu=(\mu_1,\mu_2)$, the Taylor expansions of $L(\mu)$ and ${D}(\mu)$ are
\begin{equation*}
D(\mu)=D_0+\mu_1 D_1^{(1,0)}+\mu_2 D_1^{(0,1)}+\frac{1}{2}\left(\mu_1^2 D_2^{(2,0)}+2\mu_1\mu_2 D_2^{(1,1)}+\mu_2^2 D_2^{(0,2)}\right)+\cdots,
\end{equation*}
and
\begin{equation*}
L(\mu)=L_0+\mu_1 L_1^{(1,0)}+\mu_2 L_1^{(0,1)}+\frac{1}{2}\left(\mu_1^2 L_2^{(2,0)}+2\mu_1\mu_2 L_2^{(1,1)}+\mu_2^2 L_2^{(0,2)}\right)+\cdots.
\end{equation*}
Separating the linear part from system (\ref{general reaction-diffusion equations r theta}) yields
\begin{equation}\label{Linearizing system}
\frac{\partial U(t)}{\partial t}= \tilde{L} U_t+\tilde{F}(U_t,\mu),
\end{equation}
where $\tilde{L}U=D_0\Delta_{r\theta}U+L_0U$ with $D_0=D(0),~L_0=L(0)$, and $\tilde{F}(U_t,\mu)=[D(\mu)-D_0]\Delta_{r\theta}U+[L(\mu)-L_0]U_t+F(U_t,\mu)$.

The characteristic equation of the linearized equation at zero solution of (\ref{Linearizing system}) is
\begin{equation}\label{characteristic equation}
\prod_{m_1}{\Gamma_{0m_1}(\gamma)}\prod_{n,m_2}{\tilde{\Gamma}_{nm_2}(\gamma)}=0,
\end{equation}
with
\begin{equation}\label{characteristic equations}
\begin{aligned}
&\Gamma_{0m_1}(\gamma)=\mathrm{det}\left[\gamma I +\lambda_{0m_1}D_0-L_0(\mathrm{e}^{\gamma\cdot}I)\right]=0,~m_1=0,1,2,\cdots,\\
&\tilde{\Gamma}_{nm_2}(\gamma)=\mathrm{det}\left[\gamma I +\lambda_{nm_2}D_0-L_0(\mathrm{e}^{\gamma\cdot}I)\right]^2=0,~n=1,2,\cdots,~m_2=1,2,\cdots,
\end{aligned}
\end{equation}
and
$$
  \lambda_{nm}=\left\{\begin{array}{ll}
  \frac{\alpha_{0m_1}^2}{R^2}, & m_1=0,1,2,\cdots,\\
  \frac{\alpha_{nm_2}^2}{R^2}, & n=1,2,\cdots,~m_2=1,2,\cdots,
  \end{array}
  \right.
$$
where $-\alpha_{0m_1}$ and $-\alpha_{nm_2}$ are eigenvalues of the Laplacian on the unit disk, see \citep{Murray2001M,Pinchover2005M,Chen2023J} and the corresponding unit eigenfuncitons of the Laplacian are
$$
  \hat{\phi}_{nm}^{\lambda}=\left\{\begin{array}{ll}
  \hat{\phi}_{0m_1}^{c}, & m_1=0,1,2,\cdots,\\
  \hat{\phi}_{nm_2}^{c},~\phi_{nm_2}^{s}, & n=1,2,\cdots,~m_2=1,2,\cdots,
  \end{array}
  \right.
$$
with
$$
\hat{\phi}_{0m_1}^{c}=\frac{J_{0}\left(\frac{\alpha_{0 m_1}}{R} r\right)}{\|J_{0}\left(\frac{\alpha_{0 m_1}}{R} r\right)\|},
~\hat{\phi}_{nm_2}^{c}=\frac{J_{n}\left(\frac{\alpha_{n m_2}}{R} r\right) \mathrm{e}^{\mathrm{i} n \theta}}{2\pi\|J_{n}\left(\frac{\alpha_{n m_2}}{R} r\right)\|},
~\hat{\phi}_{nm_2}^{s}=\overline{\hat{\phi}_{nm_2}^{c}}=\frac{J_{n}\left(\frac{\alpha_{n m_2}}{R} r\right) \mathrm{e}^{-\mathrm{i} n \theta}}{2\pi\|J_{n}\left(\frac{\alpha_{n m_2}}{R} r\right)\|},
$$
which form an orthonormal basis for $\mathscr{X}_\mathbb{C}$.

In order to consider the interaction of Turing instability and Hopf bifurcation, we list the following assumptions for $\mu=(0,0)$ in Table \ref{assumptionsT}.
\begin{table}
\caption{Roots with zero real part of (\ref{characteristic equation}) and the dimension of the central subspace (dim). }\label{assumptionsT}
\begin{threeparttable}
\begin{ruledtabular}
\begin{tabular}{ccccccc}
                     &  (ET-H)                       & (T-EH)                                     & (ET-EH)                               \\
   \hline
    $\Gamma_{0m_1}$ &  $\pm \mathrm{i}\omega_{H_1}$ &  0                                       &     \rule{1cm}{0.4pt}                           \\
   $\tilde{\Gamma}_{nm_2}$ &   0 (repeated)                 & $\pm \mathrm{i}\omega_{H_2}$ (repeated)    &    0 (repeated),  $\pm \mathrm{i}\omega_{H_3}$ (repeated)\\
   $m_1,n,m_2$ & $m_1=m_{H_1},n=n_{T_1},m_2=m_{T_1}$ &  $ m_1=m_{T_2},n=n_{H_2},m_2=m_{H_2}$   &   $n=n_{T_3},m_2=m_{T_3},n=n_{H_3},m_2=m_{H_3}$\\
   dim & 4&5 & 6
\end{tabular}
\begin{tablenotes}    
        \item[1] In (ET-EH), for example, the chosen indexes mean that $\tilde{\Gamma}_{n_{T_3}m_{T_3}}(0)=0,~\tilde{\Gamma}_{n_{H_3}m_{H_3}}(\pm \mathrm{i}\omega_{H_3})=0$.      
      \end{tablenotes}            
\end{ruledtabular}
\end{threeparttable}
\end{table}
Inspired by \citep{Golubitsky1989M,Guo2013M}, if (ET-H) holds, we call this is a \textit{ET-H bifurcation}, which means, the center space is spanned by the eigenvectors of a repeated zero eigenvalue (both geometric multiplicity and algebraic multiplicity are two) and a pair of simple imaginary roots.
Similarly, if (T-EH) holds, we call this a \textit{T-EH bifurcation}. If (ET-EH) holds, we call this a \textit{ET-EH bifurcation}.

\section{Main Results}\label{sec3-Normal}

In this section, based on the Turing-Hopf normal forms theory for reaction-diffusion systems in a one-dimensional interval \citep{Song2016J,Jiang2020J}, we will derive the normal forms for {ET-H}, {T-EH}, and {ET-EH bifurcations} on a disk, respectively. The normal forms for ET-H and T-EH bifurcations can be considered as parts of the normal form of the ET-EH bifurcation, and the derivation is somewhat simpler. Therefore, we first provide Theorem \ref{ET-EH reduce} on normal forms for the ET-EH bifurcation ($n_{T_3} \ne 2 n_{H_3}$) and Remark \ref{nt2nh} for $n_{T_3} = 2 n_{H_3}$, while the other two normal forms are presented as Corollaries \ref{normal form ETH} and \ref{normal form TEH}. In addition, we provide approximate forms for the solutions restricted to the center subspace corresponding to several spatiotemporal patterns in Theorem \ref{TH}, Remarks \ref{breathing patterns}, \ref{quasi-periodic solution} and \ref{RLWSLW}.
\subsection{ET-EH bifurcation}\label{A}

If (ET-EH) holds, the center subspace of the equilibrium is six-dimensional. After coordinate transformation, the normal form on the center manifolds can be transformed into a four-dimensional real ordinary differential equations (ODEs) with ${\rho}_{H^1},~{\rho}_{H^2},~{\rho}_{T^1}$ and ${\rho}_{T^2}$ as independent variables, where ${\rho}_{H^i},i=1,2$ are variables on the eigenspace corresponding to pure imaginary roots $\pm \mathrm{i}\omega_{H_3}$ (Hopf) and ${\rho}_{T^i},i=1,2$ correspond to the zero root (Turing). When $n_{T_3} \ne 2 n_{H_3}$, the detailed derivation of the normal form is presented in Appendix \ref{thproof}, and the specific transformation can be found in (\ref{trans}). When $n_{T_3} = 2 n_{H_3}$, there will be additional terms that make the normal form more complex. Therefore, we can obtain the following results.

\begin{theorem}\label{ET-EH reduce}
When $n_{T_3} \ne 2 n_{H_3}$, the normal form truncated to the third order for the ET-EH bifurcation can be written in polar coordinates as

\begin{equation}\label{rhoETEH}
\begin{aligned}
&\dot{\rho}_{H^1}=(\epsilon_1(\mu)+c_{11}\rho_{H^1}^2+c_{12}\rho_{H^2}^2+c_{13}\rho_{T^1}\rho_{T^2})\rho_{H^1},\\
&\dot{\rho}_{H^2}=(\epsilon_1(\mu)+c_{11}\rho_{H^2}^2+c_{12}\rho_{H^1}^2+c_{13}\rho_{T^1}\rho_{T^2})\rho_{H^2},\\
&\dot{\rho}_{T^1}=(\epsilon_2(\mu)+c_{21}\rho_{H^1}^2+c_{22}\rho_{H^2}^2+c_{23}\rho_{T^1}\rho_{T^2})\rho_{T^1},\\
&\dot{\rho}_{T^2}=(\epsilon_2(\mu)+c_{21}\rho_{H^1}^2+c_{22}\rho_{H^2}^2+c_{23}\rho_{T^1}\rho_{T^2})\rho_{T^2}.
\end{aligned}
\end{equation}
\end{theorem}

\begin{theorem}\label{TH}
We are mainly concerned with the properties corresponding to the following fourteen equilibrium points of (\ref{rhoETEH}), which are separated into eight categories.\\
$($ET-EH-$\mathrm{\romannumeral1})$~ $(\rho_{H^1},\rho_{H^2},\rho_{T^1},\rho_{T^2})=(0,0,0,0)$ corresponds to the origin in the six-dimensional phase space and stands for a
\textbf{stationary solution}, which is spatially homogeneous.\\
$($ET-EH-$\mathrm{\romannumeral2})$~ $(\rho_{H^1},\rho_{H^2},\rho_{T^1},\rho_{T^2})=\left(0,0,\rho_{T^1},\rho_{T^2}\right)$ with $\rho_{T^1}\rho_{T^2}=-\epsilon_2(\mu)/c_{23}$, corresponds to a \textbf{static Turing pattern}.\\
$($ET-EH-$\mathrm{\romannumeral3})$~ $(\rho_{H^1},\rho_{H^2},\rho_{T^1},\rho_{T^2})=\left(0,\sqrt{\frac{-\epsilon_1(\mu)}{c_{11}}},0,0\right)$, for $\epsilon_1(\mu)c_{11}<0$, corresponds to a periodic solution in the subspace of $(z_2,z_3)$, which is a \textbf{rotating wave solution}. At this point, the periodic solution restricted to the center subspace has the following approximate form
$$
U(t)(r,\theta) \approx
\sum_{i=1}^n{2|p_{1i}|\sqrt{\frac{-\epsilon_1(\mu)}{c_{11}}} J_{n_{H_3}}(\sqrt{\lambda_{n_{H_3}m_{H_3}}}r)\cos(\mathrm{Arg}(p_{1i})+\omega_{H_3} t+n_{H_3}\theta) \mathbf{e}_i},
$$
where $\mathbf{e}_i$ is the $i$th unit coordinate vector of $\mathbb{R}^n$ and $p_{1i}, 1\le i\le n$ are defined in Appendix \ref{thproof}.\\
$($ET-EH-$\mathrm{\romannumeral4})$~ $(\rho_{H^1},\rho_{H^2},\rho_{T^1},\rho_{T^2})=\left(\sqrt{\frac{-\epsilon_1(\mu)}{c_{11}}},0,0,0\right)$, for $\epsilon_1(\mu)c_{11}<0$,  corresponds to a a periodic solution in the subspace of $(z_1,z_4)$, which is \textbf{rotating wave solution} in the opposite direction as that in $\mathrm{(\romannumeral2)}$. At this point, the periodic solution restricted to the center subspace has the following approximate form
$$
U(t)(r,\theta) \approx
\sum_{i=1}^n{2|p_{1i}|\sqrt{\frac{-\epsilon_1(\mu)}{c_{11}}} J_{n_{H_3}}(\sqrt{\lambda_{n_{H_3}m_{H_3}}}r)\cos(\mathrm{Arg}(p_{1i})+\omega_{H_3} t-n_{H_3}\theta) {e}_i}.
$$\\
$($ET-EH-$\mathrm{\romannumeral5})$~ $(\rho_{H^1},\rho_{H^2},\rho_{T^1},\rho_{T^2})=\left(\sqrt{\frac{-\epsilon_1(\mu)}{c_{11}+c_{12}}},\sqrt{\frac{-\epsilon_1(\mu)}{c_{11}+c_{12}}},0,0\right)$ corresponds to a periodic solution, which is a \textbf{standing wave}. At this point, the periodic solution restricted to the center subspace has the following approximate form
$$
U(t)(r,\theta) \approx
\sum_{i=1}^n{4|p_{1i}|\sqrt{\frac{-\epsilon_1(\mu)}{c_{11}+c_{12}}} J_{n_{H_3}}(\sqrt{\lambda_{n_{H_3}m_{H_3}}}r)\cos(\mathrm{Arg}(p_{1i})+\omega_{H_3} t)\cos(n_{H_3}\theta) {e}_i}.
$$\\
$($ET-EH-$\mathrm{\romannumeral6})$~ $(\rho_{H^1},\rho_{H^2},\rho_{T^1},\rho_{T^2})=\left(0,\sqrt{\frac{c_{13}\epsilon_2(\mu)-c_{23}\epsilon_1(\mu)}{c_{23}c_{11}-c_{13}c_{22}}},\rho_{T^1},\rho_{T^2}\right)$
with $\rho_{T^1}\rho_{T^2}=-\frac{\epsilon_1(\mu)+c_{11}\rho_{H^2}^2}{c_{13}}$, or $\left(0,\sqrt{\frac{-\epsilon_1(\mu)}{c_{11}}},0,\rho_{T^2}\right)$ and $\left(0,\sqrt{\frac{-\epsilon_1(\mu)}{c_{11}}},\rho_{T^1},0\right)$
with $\frac{\epsilon_1(\mu)}{c_{11}}=\frac{\epsilon_2(\mu)}{c_{22}}$,
correspond to three groups of \textbf{ET-EH patterns}. At these points, the solution of real form restricted to the center subspace has the following approximate form
\begin{equation}\label{ET-EH-rotating1}
\begin{aligned}
U(t)(r,\theta) \approx
&\sum_{i=1}^n{2|p_{1i}|\rho_{H^2} J_{n_{H_3}}(\sqrt{\lambda_{n_{H_3}m_{H_3}}}r)\cos(\mathrm{Arg}(p_{1i})+\omega_{H_3} t+n_{H_3}\theta) {e}_i}\\
&+\xi_T \left(\rho_{T^1}+\rho_{T^2}\right) J_{n_{T_3}}(\sqrt{\lambda_{n_{T_3}m_{T_3}}}r)\cos(n_{T_3}\theta).\\
\end{aligned}
\end{equation}
$($ET-EH-$\mathrm{\romannumeral7})$~ $(\rho_{H^1},\rho_{H^2},\rho_{T^1},\rho_{T^2})=\left(\sqrt{\frac{c_{13}\epsilon_2(\mu)-c_{23}\epsilon_1(\mu)}{c_{23}c_{11}-c_{13}c_{21}}},0,\rho_{T^1},\rho_{T^2}\right)$
with $\rho_{T^1}\rho_{T^2}=-\frac{\epsilon_1(\mu)+c_{11}\rho_{H^1}^2}{c_{13}}$, or
$\left(\sqrt{\frac{-\epsilon_1(\mu)}{c_{11}}},0,0,\rho_{T^2}\right)$ and $\left(\sqrt{\frac{-\epsilon_1(\mu)}{c_{11}}},0,\rho_{T^1},0\right)$
with $\frac{\epsilon_1(\mu)}{c_{11}}=\frac{\epsilon_2(\mu)}{c_{22}}$,
correspond to three groups of \textbf{ET-EH patterns} in the opposite direction as that in $\mathrm{(\romannumeral6)}$. At these points, the solution restricted to the center subspace has the following approximate form
\begin{equation}\label{ET-EH-rotating2}
\begin{aligned}
U(t)(r,\theta) \approx
&\sum_{i=1}^n{2|p_{1i}|\rho_{H^1} J_{n_{H_3}}(\sqrt{\lambda_{n_{H_3}m_{H_3}}}r)\cos(\mathrm{Arg}(p_{1i})+\omega_{H_3} t-n_{H_3}\theta) {e}_i}\\
&+\xi_T \left(\rho_{T^1}+\rho_{T^2}\right) J_{n_{T_3}}(\sqrt{\lambda_{n_{T_3}m_{T_3}}}r)\cos(n_{T_3}\theta).\\
\end{aligned}
\end{equation}
$($ET-EH-$\mathrm{\romannumeral8})$~ $(\rho_{H^1},\rho_{H^2},\rho_{T^1},\rho_{T^2})=\left(\sqrt{\frac{c_{13}\epsilon_2(\mu)-c_{23}\epsilon_1(\mu)}{c_{23}(c_{11}+c_{12})-c_{13}(c_{21}+c_{22})}},\sqrt{\frac{c_{13}\epsilon_2(\mu)-c_{23}\epsilon_1(\mu)}{c_{23}(c_{11}+c_{12})-c_{13}(c_{21}+c_{22})}},\rho_{T^1},\rho_{T^2}\right)$
with $\rho_{T^1}\rho_{T^2}=-\frac{\epsilon_1(\mu)+(c_{11}+c_{12})\rho_{H^1}^2}{c_{13}}$, or
$\left(\sqrt{\frac{-\epsilon_1(\mu)}{c_{11}+c_{12}}},\sqrt{\frac{-\epsilon_1(\mu)}{c_{11}+c_{12}}},\rho_{T^1},0\right)$ and $\left(\sqrt{\frac{-\epsilon_1(\mu)}{c_{11}+c_{12}}},\sqrt{\frac{-\epsilon_1(\mu)}{c_{11}+c_{12}}},0,\rho_{T^2}\right)$ with $\frac{\epsilon_1(\mu)}{c_{11}+c_{12}}=\frac{\epsilon_2(\mu)}{c_{21}+c_{22}}$
correspond to three groups of \textbf{ET-EH patterns}. At these point, the solution restricted to the center subspace has the following approximate form
\begin{equation}\label{ET-EH-standing}
\begin{aligned}
U(t)(r,\theta) \approx
&\sum_{i=1}^n{4|p_{1i}|\rho_{H^1} J_{n_{H_3}}(\sqrt{\lambda_{n_{H_3}m_{H_3}}}r)\cos(\mathrm{Arg}(p_{1i})+\omega_{H_3} t)\cos(n_{H_3}\theta) {e}_i}\\
&+\xi_T \left(\rho_{T^1}+\rho_{T^2}\right) J_{n_{T_3}}(\sqrt{\lambda_{n_{T_3}m_{T_3}}}r)\cos(n_{T_3}\theta).\\
\end{aligned}
\end{equation}
\end{theorem}

\begin{remark}\label{ETEH-REMARK}
$($ET-EH-$\mathrm{\romannumeral6})$-$($ET-EH-$\mathrm{\romannumeral8})$ show three types of complex ET-EH patterns. We draw a schematic diagram in Figure \ref{ET-EH-rotating} of the solution in $($ET-EH-$\mathrm{\romannumeral6})$ with $n_{T_3}=1,m_{T_3}=1;~n_{H_3}=2,m_{H_3}=2$ and $\omega_{H_3}=1$ as an example, which is
\begin{equation}\label{ET-EH-rotating-example}
U(t)(r,\theta) \approx
J_{2}(\sqrt{\lambda_{22}}r)\cos(t+2\theta)+J_{1}(\sqrt{\lambda_{11}}r)\cos{\theta}.
\end{equation}
The subfigures in the first row provide ET-EH patterns like (\ref{ET-EH-rotating1}) at $t=0,~T/3,~2T/3$, and $T$, respectively, where $T\approx6$ is the period. Fixing $r=R$ and $r=R/2$, we find that despite (\ref{ET-EH-rotating-example}) is a sum of two regular patterns generating from Hopf bifurcation and  Turing instability, under the interaction of the two, the spatial form of (\ref{ET-EH-rotating-example}) is quite complex, making it difficult to summarize general rule. Similarly, the solutions in (\ref{ET-EH-rotating2}) and (\ref{ET-EH-standing}) can be explained in the same way.
\begin{figure}[htbp]
\centering
\includegraphics[width=1\textwidth]{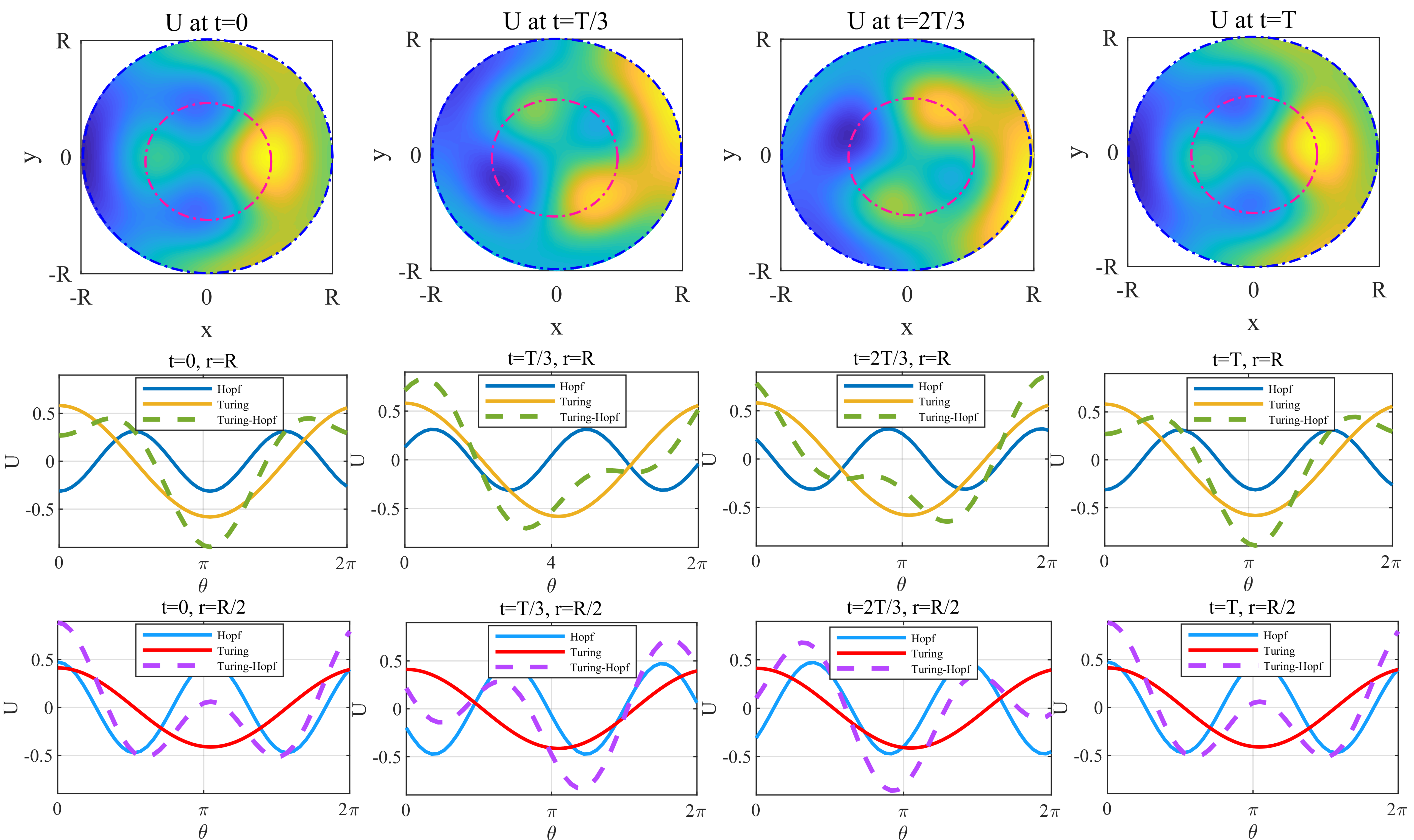}
\caption{First row: ET-EH patterns in (\ref{ET-EH-rotating-example}). Second/Third row: At $r=R/r=\frac{R}{2}$, the Hopf component, Turing component and there sum of (\ref{ET-EH-rotating-example}) are illustrated.}
\label{ET-EH-rotating}
\end{figure}
\end{remark}

\begin{remark}\label{nt2nh}
When $n_{T_3}= 2n_{H_3}$,  there will be additional terms $B_{001010}z_3z_5,~\overline{B_{001010}}z_4z_5$, $B_{001010}z_1z_6$, and $\overline{B_{001010}}z_2z_6$ in the normal form truncated to the third order for the ET-EH bifurcation. If we use the same coordinate transformation as that for $n_{T_3}\ne 2n_{H_3}$, there will be $\Delta\chi=\chi_{H^1}-\chi_{H^2}$ as a new variable, due to the presence of these additional terms. This means, by transformations $z_1=\rho_{H}\mathrm{e}^{\mathrm{i}\chi_{H^1}},~z_4=\rho_{H}\mathrm{e}^{-\mathrm{i}\chi_{H^1}},
z_3=\rho_{H}\mathrm{e}^{\mathrm{i}\chi_{H^2}},~z_2=\rho_{H}\mathrm{e}^{-\mathrm{i}\chi_{H^2}},
z_5=\rho_{T},~z_6=\rho_{T}$,
we get the normal form written in polar coordinates as
\begin{equation}\label{rhoETEH2}
\begin{aligned}
&\dot{\rho}_{H}=(\epsilon_1(\mu)+(c_{11}+c_{12})\rho_{H}^2+c_{13}\rho_{T}^2+c_{14}\rho_{T}\cos{\Delta\chi})\rho_{H},\\
&\dot{\Delta\chi}=-2c_{14}\rho_{T}\sin{\Delta\chi},\\
&\dot{\rho}_{T}=(\epsilon_2(\mu)+(c_{21}+c_{22})\rho_{H}^2+c_{23}\rho_{T}^2)\rho_{T},
\end{aligned}
\end{equation}
where $c_{14}=\mathrm{Re}\{B_{001010}\}$.

We are more concerned about the form of the original system solution corresponds to the equilibrium point of (\ref{rhoETEH2}) with ${\rho}_{H} \ne 0,~{\rho}_{T} \ne 0$ and $\Delta\chi \ne 0$, for instance, $({\rho}_{H},{\rho}_{T},\Delta\chi)=\left(\sqrt{\frac{c_{23}\rho_T^2+\epsilon_2(\mu)}{c_{21}+c_{22}}},\right.$
$\left.\frac{-C_2 \pm \sqrt{C_2^2-4C_1C_3}}{2C_1},\pi\right)$, with $C_1=(c_{11}+c_{12})c_{23}+(c_{21}+c_{22})c_{13},~C_2=\frac{-(c_{21}+c_{22})c_{14}\pi}{2}$, $C_{3}=(c_{11}+c_{12})\epsilon_2(\mu)+(c_{21}+c_{22})\epsilon_1(\mu)$.
At these points, the solution restricted to the center subspace has the following approximate form
\begin{equation}\label{ET-EH-standing2}
\begin{aligned}
U(t)(r,\theta) \approx
&-\sum_{i=1}^n{4|p_{1i}|\rho_{H} J_{n_{H_3}}(\sqrt{\lambda_{n_{H_3}m_{H_3}}}r)\sin(\mathrm{Arg}(p_{1i})+\chi_{H^1}(t))\sin(n_{H_3}\theta) {e}_i}\\
&+2\xi_T \rho_{T} J_{n_{T_3}}(\sqrt{\lambda_{n_{T_3}m_{T_3}}}r)\cos(n_{T_3}\theta).\\
\end{aligned}
\end{equation}
It can be observed that due to $\Delta\chi=\chi_{H^1}-\chi_{H^2}=\pi$, there is also a certain phase difference in the Hopf part and Turing part of the expression (\ref{ET-EH-standing2}). Thus, the form of solution maintains standing wave characteristics (Hopf) and static pattern characteristics (Turing) at different positions on the disk.
\end{remark}

\subsection{ET-H bifurcation}\label{B}

If (ET-H) holds, compared to subsection \ref{A}, the dimension of the eigenspace corresponding to pure imaginary roots $\pm\mathrm{i}\omega_{H_2}$ decreases. By (\ref{Mj1z}) and (\ref{normal forms for Turing-Hopf}), we can obtain the following results.

\begin{corollary}\label{normal form ETH}
The normal form truncated to the third order for ET-H bifurcation in polar coordinates is
\begin{equation}\label{ET-H-rho}
\begin{aligned}
&\dot{\rho}_{H}=(\alpha_1(\mu)+a_{11}\rho_{H}^2+a_{12}\rho_{T^1}\rho_{T^2})\rho_{H},\\
&\dot{\rho}_{T^1}=(\alpha_2(\mu)+a_{21}\rho_{T^1}\rho_{T^2}+a_{22}\rho_{H}^2)\rho_{T^1},\\
&\dot{\rho}_{T^2}=(\alpha_2(\mu)+a_{21}\rho_{T^1}\rho_{T^2}+a_{22}\rho_{H}^2)\rho_{T^2}.\\
\end{aligned}
\end{equation}
We can explain dynamics of the system by analyzing five equilibrium points of system (\ref{ET-H-rho}). The equilibrium points $(0,0,0)$ and $(0,\rho_{T^1},\rho_{T^2})$ with $\rho_{T^1}\rho_{T^2}=-\frac{\alpha_2(\mu)}{a_{21}}$ are similar to $($ET-EH-$\mathrm{\romannumeral1})$-$($ET-EH-$\mathrm{\romannumeral2})$, but the dynamic properties of the other equilibrium points are simpler than $($ET-EH-$\mathrm{\romannumeral6})$-$($ET-EH-$\mathrm{\romannumeral8})$. Therefore, we only provide the following remark.
\end{corollary}

\begin{remark}\label{breathing patterns}
$(\rho_{H},\rho_{T^1},\rho_{T^2})=\left(\sqrt{\frac{a_{12}\alpha_2(\mu)-a_{21}\alpha_1(\mu)}{a_{11}a_{21}-a_{12}a_{22}}},\rho_{T^1},\rho_{T^2}\right)$ with $\rho_{T^1}\rho_{T^2}=\frac{\alpha_1(\mu)+a_{11}\rho_{H}^2}{-a_{12}}$, or $(\sqrt{\frac{-\alpha_1(\mu)}{a_{11}}},\rho_{T^1},0)$ and $(\sqrt{\frac{-\alpha_1(\mu)}{a_{11}}},0,\rho_{T^1})$ with $\frac{\alpha_1(\mu)}{a_{11}}=\frac{\alpha_2(\mu)}{a_{22}}$, correspond to three groups of dynamic Turing-Hopf patterns \textbf{(breathing patterns)}. At these points, the solution restricted to the center subspace has the following approximate form
$$
\begin{aligned}
U(t)(r,\theta) \approx
&\sum_{i=1}^n{2|p_{1i}| \rho_{H} J_{0}(\sqrt{\lambda_{0m_{T_1}}}r)\cos(\mathrm{Arg}(p_{1i})+\omega_{H_1} t){e}_i}\\
&+\xi_T \left(\rho_{T^1}+\rho_{T^2}\right) J_{n_{T_1}}(\sqrt{\lambda_{n_{T_1}m_{T_1}}}r)\cos(n_{T_1}\theta).
\end{aligned}
$$
Similar to the discussion in Remark \ref{ETEH-REMARK}, the solution will maintain a fixed inhomogeneous form and oscillate up and down over time (breathing).
\end{remark}

\begin{remark}
Let $\rho_{T}^2=\rho_{T^1}\rho_{T^2},~\bar{\rho}_{H}={\rho}_{H}\sqrt{|a_{11}|},~~\bar{\rho}_{T}={\rho}_{T}\sqrt{|a_{21}|}$,  and drop the bars, then system (\ref{ET-H-rho}) can be transformed into
\begin{equation}\label{rhoHrhoT}
\begin{aligned}
&\dot{\rho}_{H}=(\alpha_1(\mu)+\rho_{H}^2+a_b\rho_{T}^2)\rho_{H},\\
&\dot{\rho}_{T}=(\alpha_2(\mu)+a_c\rho_{H}^2+a_d\rho_{T}^2)\rho_{T},\\
\end{aligned}
\end{equation}
which has twelve distinct kinds of unfoldings. The stability conditions of equilibrium points can be given, by Chapter 7.5 in \citep{Guckenheimer1983M}. Thus, in this case, the stability of spatiotemporal solutions and a complete bifurcation set are easily obtained.
\end{remark}

\begin{remark}\label{quasi-periodic solution}
By the case \uppercase\expandafter{\romannumeral6}a of Chapter 7.5 in \citep{Guckenheimer1983M}, there is a quasi-periodic solution on the three-dimensional torus, which corresponds to that system (\ref{rhoHrhoT}) has a center and level curves with $\rho_{H}^2+\upsilon\rho_{T}^2=-\alpha_1(\mu)$ where $\upsilon=\frac{a_b+1}{a_c-1}$. The solution generated by the Hopf bifurcation restricted to the center subspace has the following approximate form
$$
\begin{aligned}
U(t)(r,\theta) &\approx \sum_{i=1}^n{2|p_{1i}| \rho_{H} J_{0}(\sqrt{\lambda_{0m_{T_1}}}r)\cos(\mathrm{Arg}(p_{1i})+\omega_{H_1} t)\cos(\bar{\omega} t){e}_i}\\
&+\xi_T \rho_{T} J_{n_{T_1}}(\sqrt{\lambda_{n_{T_1}m_{T_1}}}r)\cos(n_{T_1}\theta)\sin(\bar{\omega} t),
\end{aligned}
$$
where $\bar{\omega}=O(\alpha_i(\mu))$. This is a rather complicated pattern including one spatial frequency and two different temporal frequencies, which is actually a quasi-periodic oscillation with spatial inhomogeneous profiles.
\end{remark}

\subsection{T-EH bifurcation}

If (T-EH) holds, compared to subsection \ref{A}, the dimension of the eigenspace corresponding to the zero root decreases and the following results can be obtained.
\begin{corollary}\label{normal form TEH}
The normal form truncated to the third order for T-EH bifurcation in polar coordinates is
\begin{equation}\label{T-EH-rho}
\begin{aligned}
&\dot{\rho}_{H^1}=(\beta_1(\mu)+b_1{\rho}_{T}+b_{11}\rho_{H^1}^2+b_{12}\rho_{H^2}^2+b_{13}\rho_{T}^2)\rho_{H^1},\\
&\dot{\rho}_{H^2}=(\beta_1(\mu)+b_1{\rho}_{T}+b_{11}\rho_{H^2}^2+b_{12}\rho_{H^1}^2+b_{13}\rho_{T}^2)\rho_{H^2},\\
&\dot{\rho}_{T}=(\beta_2(\mu)+b_2{\rho}_{T}+b_{21}\rho_{H^1}^2+b_{22}\rho_{H^2}^2+b_{23}\rho_{T}^2)\rho_{T}.\\
\end{aligned}
\end{equation}
We can explain dynamics of the system by analyzing at most twelve equilibrium points of system (\ref{T-EH-rho}). Similar to subsection \ref{B}, several equilibrium points of system (\ref{T-EH-rho}) are consistent with the results of Theorem \ref{TH}. Next, we will explain in detail several solutions for the interaction of Turing-Hopf under {\rm{(T-EH)}}, which is more clearer than $($ET-EH-$\mathrm{\romannumeral6})$-$($ET-EH-$\mathrm{\romannumeral8})$.
\end{corollary}

\begin{remark}\label{RLWSLW}
$\mathrm{(\romannumeral1)}$~ $(\rho_{H^1},\rho_{H^2},\rho_{T})=\left(0,\sqrt{\frac{\beta_1(\mu)+b_1\rho_{T}+b_{13}\rho_{T}^2}{-b_{11}}},\frac{-B_2\pm\sqrt{B_2^2-4B_1B_3}}{2B_1}\right)$ with $B_1=b_1b_{22}-b_2b_{11},~B_2=b_{22}b_{13}-b_{11}b_{23},~B_3=b_{22}\beta_1(\mu)-b_{11}\beta_2(\mu)$, correspond to at most two \textbf{rotating wave-like dynamic Turing-Hopf patterns}, depending on the sign of $B_2^2-4B_1B_3$. At these points, the periodic solution restricted to the center subspace has the following approximate form
$$
\begin{aligned}
U(t)(r,\theta) \approx
&\sum_{i=1}^n{2|p_{1i}|\rho_{H^2} J_{n_{H_2}}(\sqrt{\lambda_{n_{H_2}m_{H_2}}}r)\cos(\mathrm{Arg}(p_{1i})+\omega_{H_2} t+n_{H_2}\theta) {e}_i}\\
&+\xi_T \rho_{T} J_0(\sqrt{\lambda_{0m_{T_2}}}r).
\end{aligned}
$$
Similarly, the spatial form of the Turing component is constant. Therefore, along with a circle with radius $r$ on the disk, the solution will be in the form of a rotating wave.\\
$\mathrm{(\romannumeral2)}$~ $(\rho_{H^1},\rho_{H^2},\rho_{T})=\left(\sqrt{\frac{\beta_1(\mu)+b_1\rho_{T}+b_{13}\rho_{T}^2}{-b_{11}}},0,\frac{-B_5\pm\sqrt{B_5^2-4B_4B_6}}{2B_4}\right)$ with $B_4=b_1b_{21}-b_2b_{11},~B_5=b_{21}b_{13}-b_{11}b_{23},~B_6=b_{21}\beta_1(\mu)-b_{11}\beta_2(\mu)$, correspond to at most two \textbf{rotating wave-like dynamic Turing-Hopf patterns} in the opposite direction as that in $\mathrm{(\romannumeral1)}$. At these points, the periodic solution restricted to the center subspace has the following approximate form
$$
\begin{aligned}
U(t)(r,\theta) \approx
&\sum_{i=1}^n{2|p_{1i}|\rho_{H^2} J_{n_{H_2}}(\sqrt{\lambda_{n_{H_2}m_{H_2}}}r)\cos(\mathrm{Arg}(p_{1i})+\omega_{H_2} t-n_{H_2}\theta) {e}_i}\\
&+\xi_T \rho_{T} J_0(\sqrt{\lambda_{0m_{T_2}}}r).
\end{aligned}
$$\\
$\mathrm{(\romannumeral3)}$~ $(\rho_{H^1},\rho_{H^2},\rho_{T})=\left(\sqrt{\frac{\beta_1(\mu)+b_1\rho_{T}+b_{13}\rho_{T}^2}{-(b_{11}+b_{12})}},\sqrt{\frac{\beta_1(\mu)+b_1\rho_{T}+b_{13}\rho_{T}^2}{-(b_{11}+b_{12})}},\frac{-B_8\pm\sqrt{B_8^2-4B_7B_9}}{2B_7}\right)$ with $B_7=(b_{21}+b_{22})b_{13}-(b_{11}+b_{12})b_{23},~B_8=(b_{21}+b_{22})b_{1}-(b_{11}+b_{12})b_{2},~B_9=(b_{21}+b_{22})\beta_1(\mu)-(b_{11}+b_{12})\beta_2(\mu)$, correspond to at most two \textbf{standing wave-like dynamic Turing-Hopf patterns}. At these points, the periodic solution restricted to the center subspace has the following approximate form
$$
\begin{aligned}
U(t)(r,\theta) \approx
&\sum_{i=1}^n{4|p_{1i}| \rho_{H^1} J_{n_{H_2}}(\sqrt{\lambda_{n_{H_2}m_{H_2}}}r)\cos(\mathrm{Arg}(p_{1i})+\omega_{H_2} t)\cos(n_{H_2}\theta) {e}_i}\\
&+\xi_T \rho_{T} J_0(\sqrt{\lambda_{0m_{T_2}}}r).
\end{aligned}
$$
\end{remark}


\section{Numerical Simulations}\label{sec4-Numerical Simulations}
In \citep{Shen2019J}, Shen and Wei investigated a delayed mussel-algae system. Here, we investigate the dynamics of such a model on a disk.
\begin{equation}\label{mussel-algae}
\left\{\begin{array}{l}
\frac{\partial m(t, r, \theta)}{\partial t}=d_{1} \Delta_{r \theta} m(t, r, \theta)+m(t,r,\theta)\left(ba(t-\tau,r,\theta)-\frac{1}{1-m(t-\tau,r,\theta)}\right),~(r, \theta) \in \mathbb{D},~t>0, \\
\kappa\frac{\partial a(t, r, \theta)}{\partial t}=\Delta_{r \theta} a(t, r, \theta)+\alpha\left(1-a(t,r,\theta)\right) -m(t,r,\theta)a(t,r,\theta),~(r, \theta) \in \mathbb{D},~t>0, \\
\partial_{r} m(\cdot, R, \theta)=\partial_{r} a(\cdot, R, \theta)=0,~\theta \in [0,2\pi),
\end{array}\right.
\end{equation}
where the parameters are defined in \citep{Shen2019J}. In real-world, limited source, like nutrients and light, can lead to nonlocal intraspecific competition among algae in the ocean \citep{Steen2003J,Manoylov2009J}. Therefore, based on system (\ref{mussel-algae}), we introduced nonlocal effects by replacing $\alpha\left(1-a(t,r,\theta)\right)$ by $\alpha\left(1-\hat{a}(t,r,\theta)\right)$ with
$$
\hat{a}(t,r,\theta)=\frac{1}{ \pi R^2} \int_{0}^{R} \int_{0}^{2 \pi} \bar{r} a\left(t,\bar{r},\bar{\theta}\right)\mathrm{d} \bar{\theta} \mathrm{d} \bar{r}.
$$
Then, system (\ref{mussel-algae}) becomes
\begin{equation}\label{mussel-algae with nonlocal}
\left\{\begin{array}{l}
\frac{\partial m(t, r, \theta)}{\partial t}=d_{1} \Delta_{r \theta} m(t, r, \theta)+m(t,r,\theta)\left(ba(t-\tau,r,\theta)-\frac{1}{1-m(t-\tau,r,\theta)}\right),~(r, \theta) \in \mathbb{D},~t>0, \\
\kappa\frac{\partial a(t, r, \theta)}{\partial t}=\Delta_{r \theta} a(t, r, \theta)+\alpha\left(1-\hat{a}(t,r,\theta)\right) -m(t,r,\theta)a(t,r,\theta),~(r, \theta) \in \mathbb{D},~t>0, \\
\partial_{r} m(\cdot, R, \theta)=\partial_{r} a(\cdot, R, \theta)=0,~\theta \in [0,2\pi).
\end{array}\right.
\end{equation}

\begin{figure}[htbp]
\centering
\includegraphics[width=0.92\textwidth]{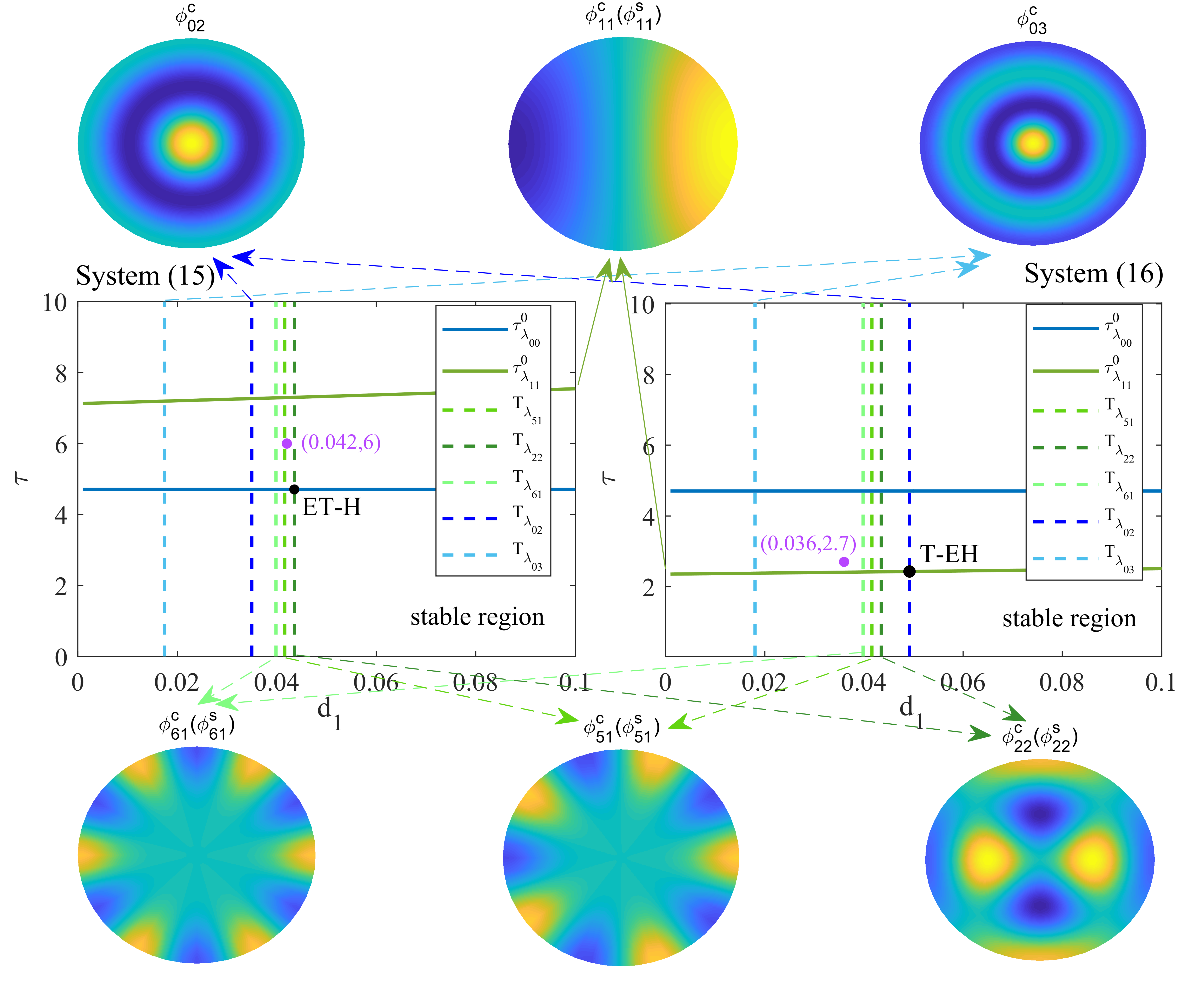}
\caption{Partial bifurcation curves on the $d_1-\tau$ plane for two systems and eigenfuncitons related to Turing instability.}
\label{bifurcation}
\end{figure}

\begin{figure}[htbp]
\centering
(a)\includegraphics[width=0.45\textwidth]{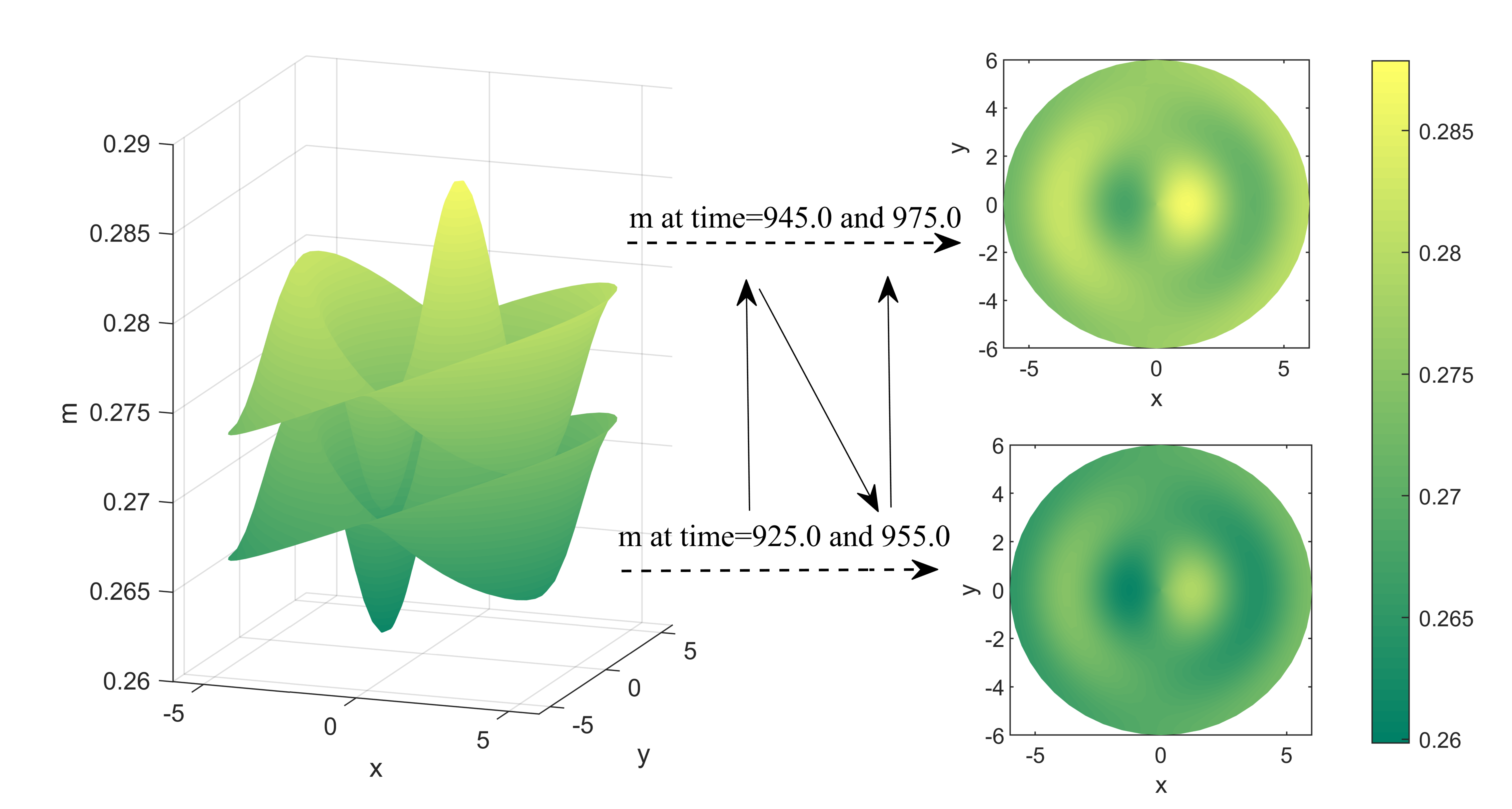}
(b)\includegraphics[width=0.45\textwidth]{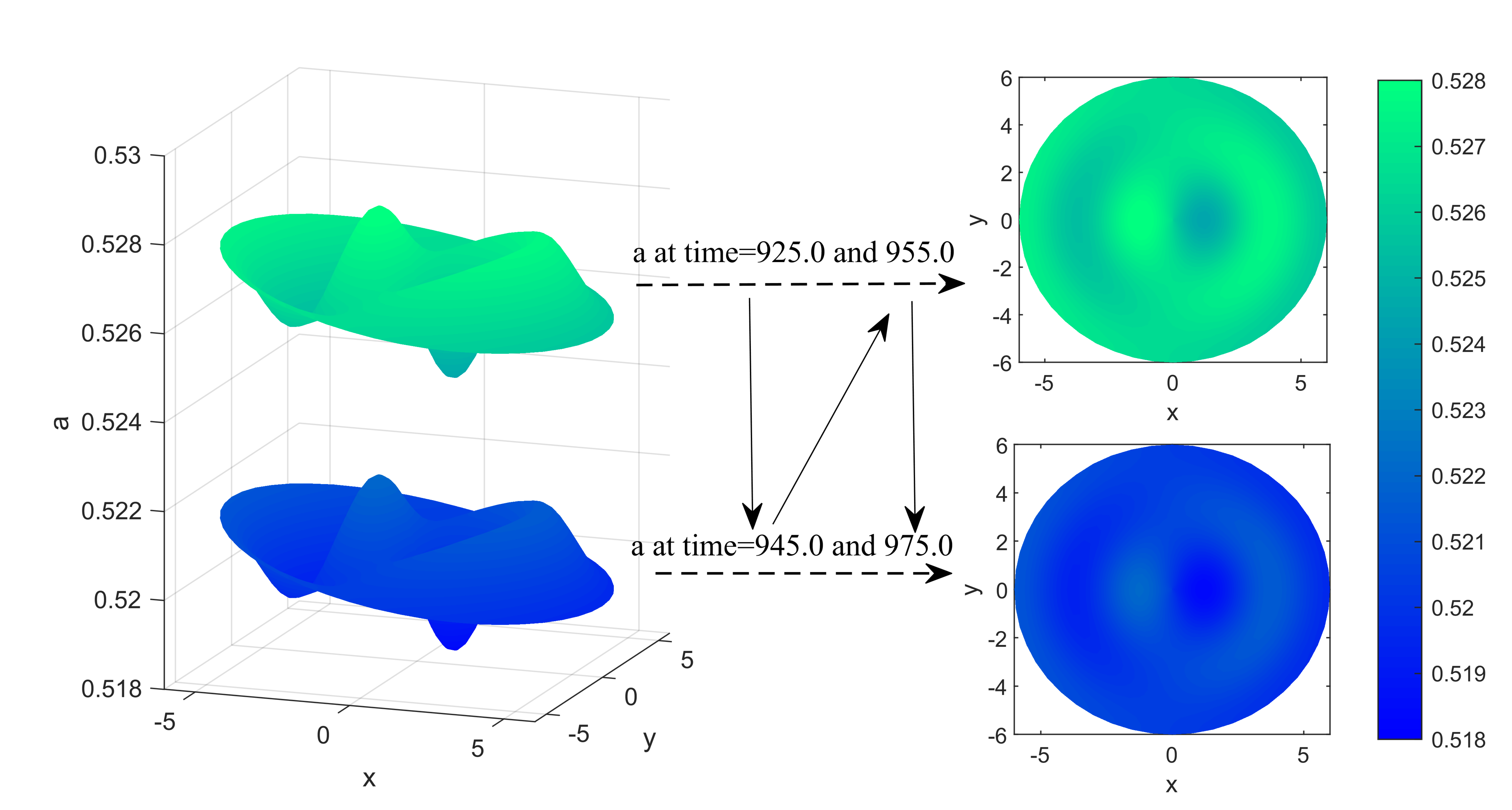}
\caption{System (\ref{mussel-algae}) produces breathing patterns with parameters: $b=1.5,~\kappa=1,~\alpha=0.3,~R=6,~d_1=0.042,~\tau=6$. Initial values are $m(t,r,\theta)=0.2727+0.01\cdot \cos t\cdot \cos r\cdot  \cos\theta,~a(t,r,\theta)=0.5238+0.01\cdot \cos t\cdot \cos r\cdot  \cos\theta,~t\in[-\tau,0)$. $(a):m,~(b):a$.}
\label{B-1}
\end{figure}

\begin{figure}[htbp]
\centering
(a)\includegraphics[width=0.45\textwidth]{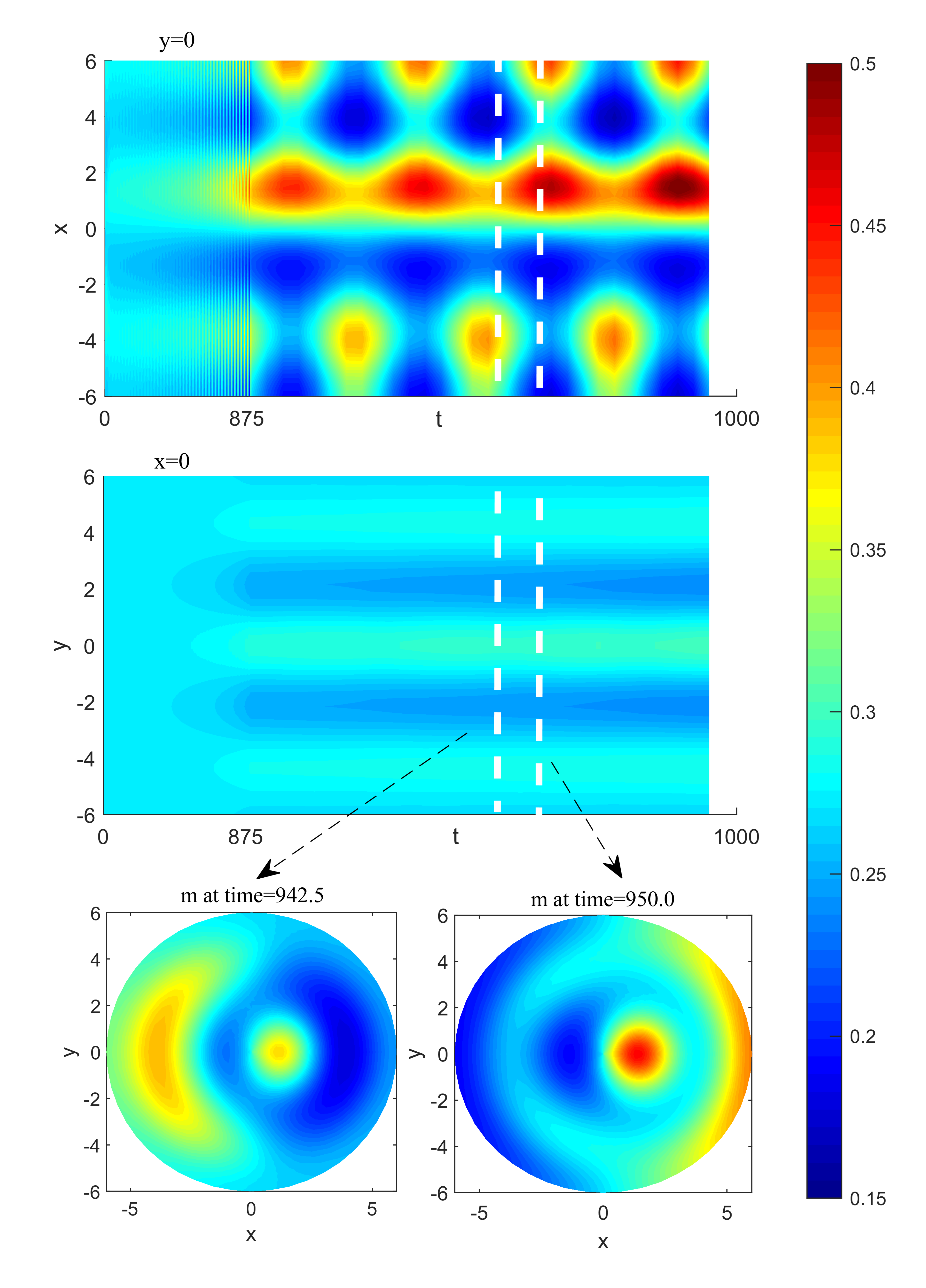}
(b)\includegraphics[width=0.45\textwidth]{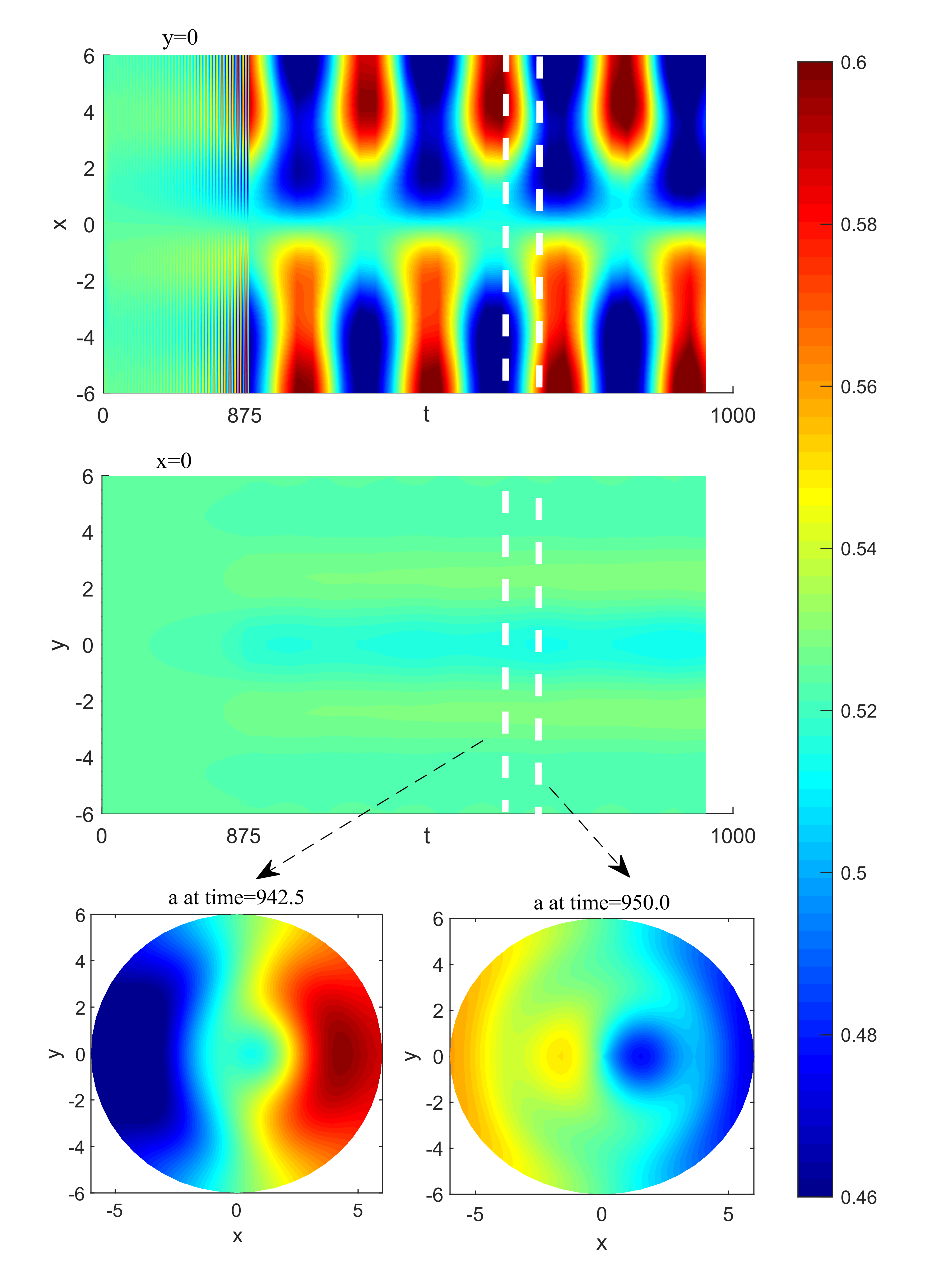}
\caption{System (\ref{mussel-algae with nonlocal}) produces standing wave-like dynamic Turing-Hopf patterns with parameters: $b=1.5,~\kappa=1,~\alpha=0.3,~R=6,~d_1=0.036,~\tau=2.7$. Initial values are $m(t,r,\theta)=0.2727+0.01\cdot \cos t\cdot \cos r\cdot  \cos\theta,~a(t,r,\theta)=0.5238+0.01\cdot \cos t\cdot \cos r\cdot  \cos\theta,~t\in[-\tau,0)$. $(a):m,~(b):a$.}
\label{SL-1}
\end{figure}

\begin{figure}[htbp]
\centering
(a)\includegraphics[width=0.45\textwidth]{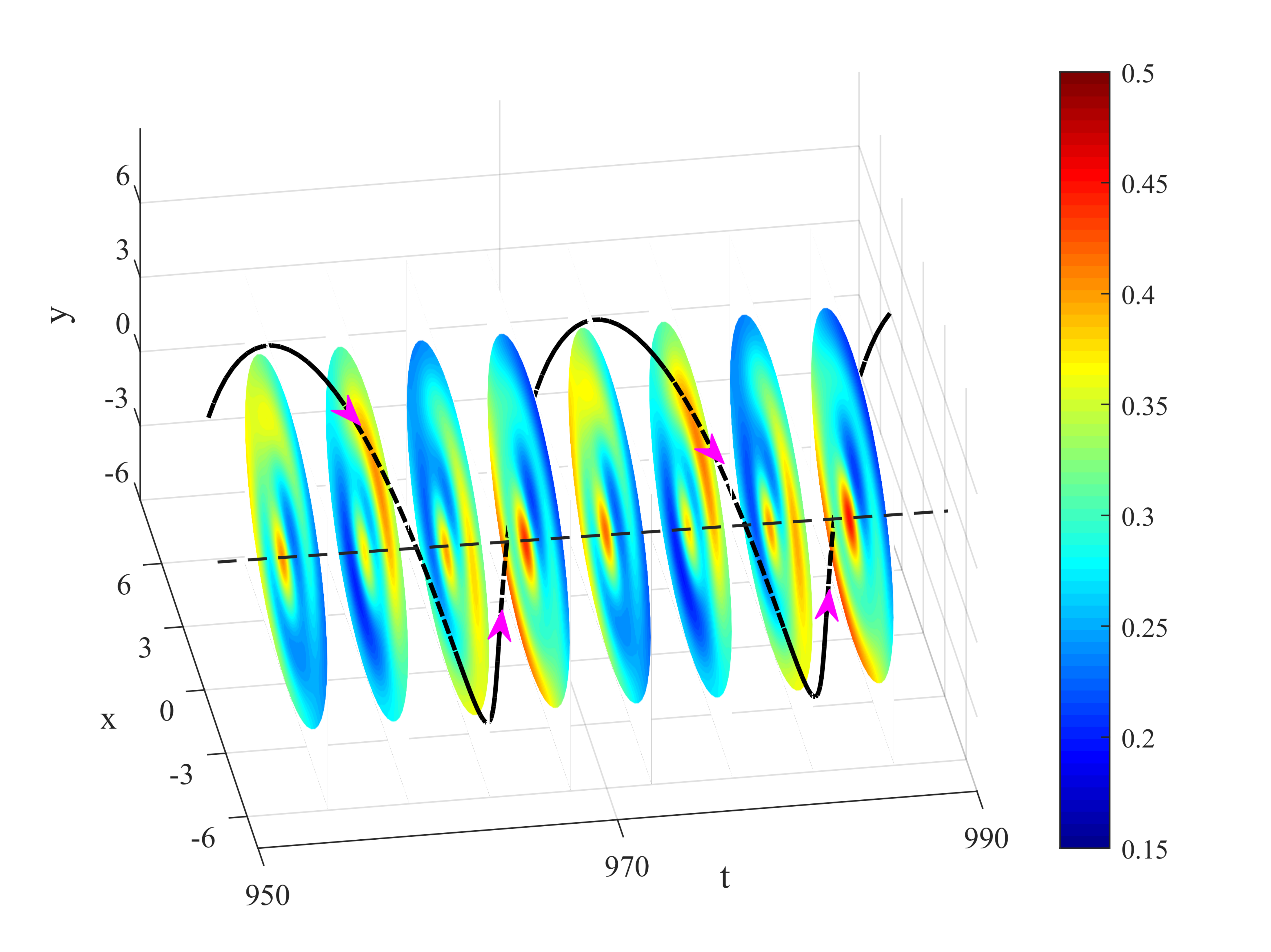}
(b)\includegraphics[width=0.45\textwidth]{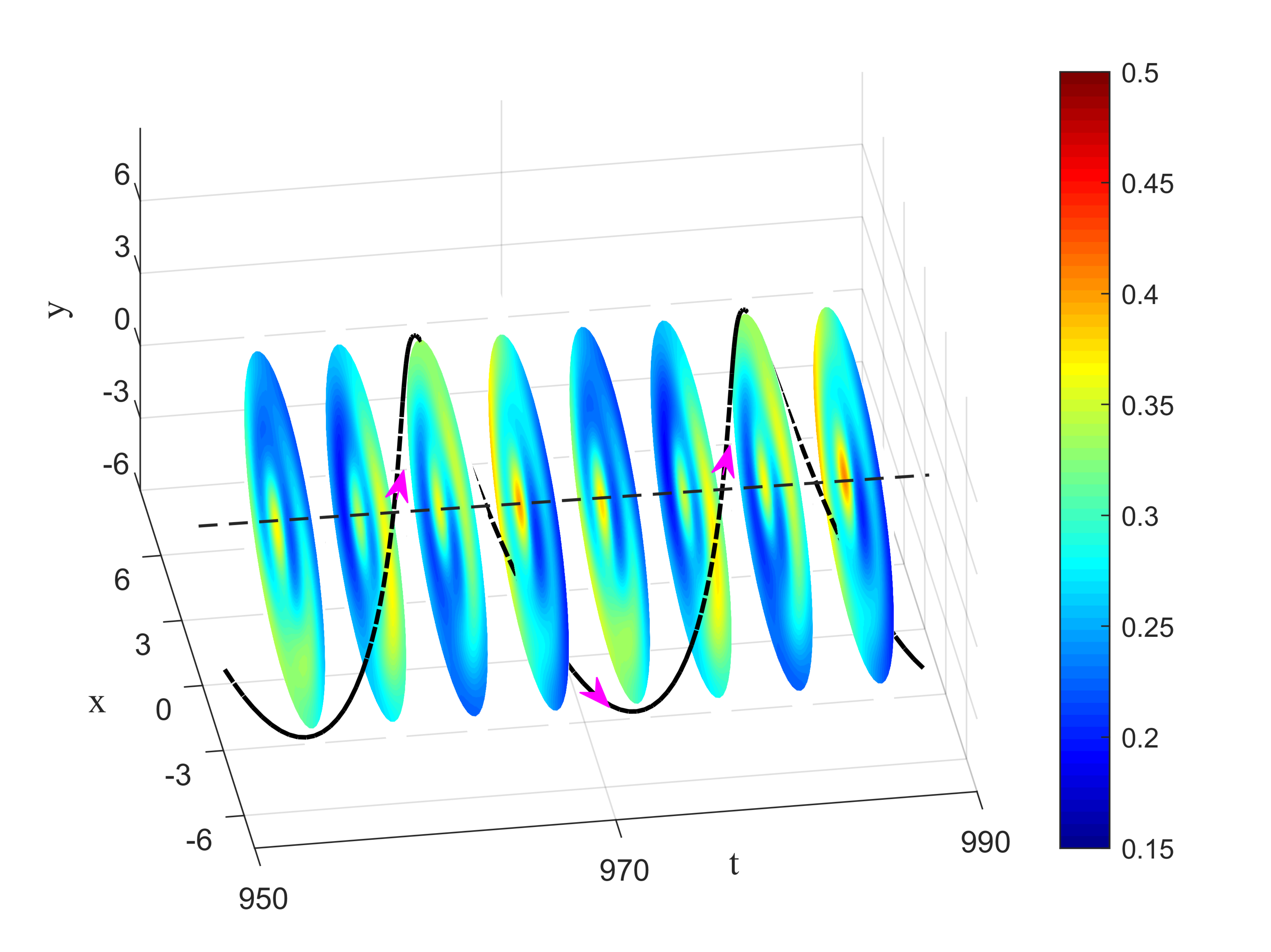}
\caption{Rotating wave-like dynamic Turing-Hopf patterns of $m$ with parameters: $b=1.5,~\kappa=1,~\alpha=0.3,~R=6,~d_1=0.036,~\tau=2.7$. Initial values are $m(t,r,\theta)=0.2727+0.01\cdot \cos t\cdot \cos r\cdot  \Theta_1(\theta),~a(t,r,\theta)=0.5238+0.01\cdot \cos t\cdot \cos r\cdot  \Theta_2(\theta),~t\in[-\tau,0)$. $(a):(\Theta_1(\theta),\Theta_2(\theta))=(\cos\theta,\sin\theta)-clockwise,~(b):(\Theta_1(\theta),\Theta_2(\theta))=(\sin\theta,\cos\theta)-anticlockwise$.}
\label{RL-1}
\end{figure}

Fixing $b=1.5,~\kappa=1,~\alpha=0.3,~R=6$, we obtain partial bifurcation curves on the $d_1-\tau$ plane of system (\ref{mussel-algae}) and system (\ref{mussel-algae with nonlocal}) shown in Figure \ref{bifurcation}, respectively. For system (\ref{mussel-algae}),
we select $(d_1,\tau)=(0.042,6)$ and get a type of breathing patterns (see Figure \ref{B-1}). For system (\ref{mussel-algae with nonlocal}),
we select $(d_1,\tau)=(0.036,2.7)$, and get two different types of dynamic Turing-Hopf patterns. Similar to the results in \citep{Chen2023J},
Turing-Hopf pattern is standing wave-like with a specific initial value (see Figure \ref{SL-1}), and with other initial values, rotating wave-like Turing-Hopf patterns appear (see Figure \ref{RL-1}).

The standing wave-like pattern has a fixed axis (see the subgraph corresponding to $y=0$ in Figure \ref{SL-1}) and a hot/cold spot indicating local maximum/minimum that does not change position over time (see the area on the right side of the fixed axis). The other parts of the pattern oscillate in the form of standing waves on both sides of the fixed axis (as shown in the subgraph corresponding to $x=0$ in Figure \ref{SL-1}).
The rotating wave-like pattern in Figure \ref{RL-1} has a portion of the pattern that remains unchanged in position and the other parts of the pattern that change in the form of rotating wave.

\section{Concluding Remarks}
In this paper, we investigate the interaction of Turing instability and Hopf bifurcation on a disk. We first present three Turing-Hopf normal forms based on different types of eigenspaces. In addition, we analyzed the possible solutions for each normal form, which can guide us to find solutions with physical significance in real-world systems. Finally, breathing, standing wave-like, and rotating wave-like patterns were simulated in a specific mussel-algae model.

Under the case (ET-EH), the possible solutions are complex, and there are several questions that can be further discussed. We believe that quasi-periodic solutions may also exist, which is quite difficult to study. In addition, in previous studies on double Hopf bifurcation, the resonance may occur: if the ratio of two imaginary roots $\mathrm{i}\omega_1$ and $\mathrm{i}\omega_2$ is rational, some additional terms cannot be eliminated. In this paper, another kind of resonance of Turing and Hopf appears, i.e. $n_{T_3}=2n_{H_3}$. Combining these factors and investigating the corresponding normal forms is a noteworthy issue to be further considered.

\bibliography{TuringHopf}

\begin{thebibliography}{40}%
\makeatletter
\providecommand \@ifxundefined [1]{%
 \@ifx{#1\undefined}
}%
\providecommand \@ifnum [1]{%
 \ifnum #1\expandafter \@firstoftwo
 \else \expandafter \@secondoftwo
 \fi
}%
\providecommand \@ifx [1]{%
 \ifx #1\expandafter \@firstoftwo
 \else \expandafter \@secondoftwo
 \fi
}%
\providecommand \natexlab [1]{#1}%
\providecommand \enquote  [1]{``#1''}%
\providecommand \bibnamefont  [1]{#1}%
\providecommand \bibfnamefont [1]{#1}%
\providecommand \citenamefont [1]{#1}%
\providecommand \href@noop [0]{\@secondoftwo}%
\providecommand \href [0]{\begingroup \@sanitize@url \@href}%
\providecommand \@href[1]{\@@startlink{#1}\@@href}%
\providecommand \@@href[1]{\endgroup#1\@@endlink}%
\providecommand \@sanitize@url [0]{\catcode `\\12\catcode `\$12\catcode
  `\&12\catcode `\#12\catcode `\^12\catcode `\_12\catcode `\%12\relax}%
\providecommand \@@startlink[1]{}%
\providecommand \@@endlink[0]{}%
\providecommand \url  [0]{\begingroup\@sanitize@url \@url }%
\providecommand \@url [1]{\endgroup\@href {#1}{\urlprefix }}%
\providecommand \urlprefix  [0]{URL }%
\providecommand \Eprint [0]{\href }%
\providecommand \doibase [0]{http://dx.doi.org/}%
\providecommand \selectlanguage [0]{\@gobble}%
\providecommand \bibinfo  [0]{\@secondoftwo}%
\providecommand \bibfield  [0]{\@secondoftwo}%
\providecommand \translation [1]{[#1]}%
\providecommand \BibitemOpen [0]{}%
\providecommand \bibitemStop [0]{}%
\providecommand \bibitemNoStop [0]{.\EOS\space}%
\providecommand \EOS [0]{\spacefactor3000\relax}%
\providecommand \BibitemShut  [1]{\csname bibitem#1\endcsname}%
\let\auto@bib@innerbib\@empty
\bibitem [{\citenamefont {Camara}\ \emph {et~al.}(2016)\citenamefont {Camara},
  \citenamefont {Haque},\ and\ \citenamefont {Mokrani}}]{Camara2016J}%
  \BibitemOpen
  \bibfield  {author} {\bibinfo {author} {\bibfnamefont {B.~I.}\ \bibnamefont
  {Camara}}, \bibinfo {author} {\bibfnamefont {M.}~\bibnamefont {Haque}}, \
  and\ \bibinfo {author} {\bibfnamefont {H.}~\bibnamefont {Mokrani}},\
  }\href@noop {} {\bibfield  {journal} {\bibinfo  {journal} {Physica A:
  Statistical Mechanics and its Applications}\ }\textbf {\bibinfo {volume}
  {461}},\ \bibinfo {pages} {374} (\bibinfo {year} {2016})}\BibitemShut
  {NoStop}%
\bibitem [{\citenamefont {Yang}\ and\ \citenamefont {Song}(2016)}]{Yang2016J}%
  \BibitemOpen
  \bibfield  {author} {\bibinfo {author} {\bibfnamefont {R.}~\bibnamefont
  {Yang}}\ and\ \bibinfo {author} {\bibfnamefont {Y.}~\bibnamefont {Song}},\
  }\href@noop {} {\bibfield  {journal} {\bibinfo  {journal} {Nonlinear
  Analysis: Real World Applications}\ }\textbf {\bibinfo {volume} {31}},\
  \bibinfo {pages} {356} (\bibinfo {year} {2016})}\BibitemShut {NoStop}%
\bibitem [{\citenamefont {Cao}\ and\ \citenamefont {Jiang}(2018)}]{Cao2018J}%
  \BibitemOpen
  \bibfield  {author} {\bibinfo {author} {\bibfnamefont {X.}~\bibnamefont
  {Cao}}\ and\ \bibinfo {author} {\bibfnamefont {W.}~\bibnamefont {Jiang}},\
  }\href@noop {} {\bibfield  {journal} {\bibinfo  {journal} {Nonlinear
  Analysis: Real World Applications}\ }\textbf {\bibinfo {volume} {43}},\
  \bibinfo {pages} {428} (\bibinfo {year} {2018})}\BibitemShut {NoStop}%
\bibitem [{\citenamefont {Kumari}\ and\ \citenamefont
  {Mohan}(2020)}]{Kumari2020J}%
  \BibitemOpen
  \bibfield  {author} {\bibinfo {author} {\bibfnamefont {N.}~\bibnamefont
  {Kumari}}\ and\ \bibinfo {author} {\bibfnamefont {N.}~\bibnamefont {Mohan}},\
  }\href@noop {} {\bibfield  {journal} {\bibinfo  {journal} {Nonlinear
  Dynamics}\ }\textbf {\bibinfo {volume} {100}},\ \bibinfo {pages} {763}
  (\bibinfo {year} {2020})}\BibitemShut {NoStop}%
\bibitem [{\citenamefont {Kumar}\ and\ \citenamefont
  {Gangopadhyay}(2020)}]{Kumar2020J}%
  \BibitemOpen
  \bibfield  {author} {\bibinfo {author} {\bibfnamefont {P.}~\bibnamefont
  {Kumar}}\ and\ \bibinfo {author} {\bibfnamefont {G.}~\bibnamefont
  {Gangopadhyay}},\ }\href@noop {} {\bibfield  {journal} {\bibinfo  {journal}
  {Physical Review E}\ }\textbf {\bibinfo {volume} {101}},\ \bibinfo {pages}
  {042204} (\bibinfo {year} {2020})}\BibitemShut {NoStop}%
\bibitem [{\citenamefont {Cross}\ and\ \citenamefont
  {Hohenberg}(1993)}]{Cross1993J}%
  \BibitemOpen
  \bibfield  {author} {\bibinfo {author} {\bibfnamefont {M.~C.}\ \bibnamefont
  {Cross}}\ and\ \bibinfo {author} {\bibfnamefont {P.~C.}\ \bibnamefont
  {Hohenberg}},\ }\href@noop {} {\bibfield  {journal} {\bibinfo  {journal}
  {Review of Modern Physics}\ }\textbf {\bibinfo {volume} {65}},\ \bibinfo
  {pages} {851} (\bibinfo {year} {1993})}\BibitemShut {NoStop}%
\bibitem [{\citenamefont {Perraud}\ \emph {et~al.}(1993)\citenamefont
  {Perraud}, \citenamefont {De~Wit}, \citenamefont {Dulos}, \citenamefont
  {De~Kepper}, \citenamefont {Dewel},\ and\ \citenamefont
  {Borckmans}}]{Perraud1993J}%
  \BibitemOpen
  \bibfield  {author} {\bibinfo {author} {\bibfnamefont {J.~J.}\ \bibnamefont
  {Perraud}}, \bibinfo {author} {\bibfnamefont {A.}~\bibnamefont {De~Wit}},
  \bibinfo {author} {\bibfnamefont {E.}~\bibnamefont {Dulos}}, \bibinfo
  {author} {\bibfnamefont {P.}~\bibnamefont {De~Kepper}}, \bibinfo {author}
  {\bibfnamefont {G.}~\bibnamefont {Dewel}}, \ and\ \bibinfo {author}
  {\bibfnamefont {P.}~\bibnamefont {Borckmans}},\ }\href@noop {} {\bibfield
  {journal} {\bibinfo  {journal} {Physical Review Letters}\ }\textbf {\bibinfo
  {volume} {71}},\ \bibinfo {pages} {1272} (\bibinfo {year}
  {1993})}\BibitemShut {NoStop}%
\bibitem [{\citenamefont {Heidemann}\ \emph {et~al.}(1993)\citenamefont
  {Heidemann}, \citenamefont {Bode},\ and\ \citenamefont
  {Purwins}}]{Heidemann1993J}%
  \BibitemOpen
  \bibfield  {author} {\bibinfo {author} {\bibfnamefont {G.}~\bibnamefont
  {Heidemann}}, \bibinfo {author} {\bibfnamefont {M.}~\bibnamefont {Bode}}, \
  and\ \bibinfo {author} {\bibfnamefont {H.~G.}\ \bibnamefont {Purwins}},\
  }\href@noop {} {\bibfield  {journal} {\bibinfo  {journal} {Physics Letters
  A}\ }\textbf {\bibinfo {volume} {177}},\ \bibinfo {pages} {225} (\bibinfo
  {year} {1993})}\BibitemShut {NoStop}%
\bibitem [{\citenamefont {Vallette}\ \emph {et~al.}(1994)\citenamefont
  {Vallette}, \citenamefont {Edwards},\ and\ \citenamefont
  {Gollub}}]{Vallette1994J}%
  \BibitemOpen
  \bibfield  {author} {\bibinfo {author} {\bibfnamefont {D.~P.}\ \bibnamefont
  {Vallette}}, \bibinfo {author} {\bibfnamefont {W.~S.}\ \bibnamefont
  {Edwards}}, \ and\ \bibinfo {author} {\bibfnamefont {J.~P.}\ \bibnamefont
  {Gollub}},\ }\href@noop {} {\bibfield  {journal} {\bibinfo  {journal}
  {Physical Review E}\ }\textbf {\bibinfo {volume} {49}},\ \bibinfo {pages}
  {R4783} (\bibinfo {year} {1994})}\BibitemShut {NoStop}%
\bibitem [{\citenamefont {Song}\ and\ \citenamefont {Zou}(2014)}]{Song2014J}%
  \BibitemOpen
  \bibfield  {author} {\bibinfo {author} {\bibfnamefont {Y.}~\bibnamefont
  {Song}}\ and\ \bibinfo {author} {\bibfnamefont {X.}~\bibnamefont {Zou}},\
  }\href@noop {} {\bibfield  {journal} {\bibinfo  {journal} {Computers \&
  Mathematics with Applications}\ }\textbf {\bibinfo {volume} {67}},\ \bibinfo
  {pages} {1978} (\bibinfo {year} {2014})}\BibitemShut {NoStop}%
\bibitem [{\citenamefont {An}\ and\ \citenamefont {Jiang}(2018)}]{An2018J}%
  \BibitemOpen
  \bibfield  {author} {\bibinfo {author} {\bibfnamefont {Q.}~\bibnamefont
  {An}}\ and\ \bibinfo {author} {\bibfnamefont {W.}~\bibnamefont {Jiang}},\
  }\href@noop {} {\bibfield  {journal} {\bibinfo  {journal} {International
  Journal of Bifurcation and Chaos}\ }\textbf {\bibinfo {volume} {28}},\
  \bibinfo {pages} {1850108} (\bibinfo {year} {2018})}\BibitemShut {NoStop}%
\bibitem [{\citenamefont {Rovinsky}\ and\ \citenamefont
  {Menzinger}(1992)}]{Rovinsky1992J}%
  \BibitemOpen
  \bibfield  {author} {\bibinfo {author} {\bibfnamefont {A.}~\bibnamefont
  {Rovinsky}}\ and\ \bibinfo {author} {\bibfnamefont {M.}~\bibnamefont
  {Menzinger}},\ }\href@noop {} {\bibfield  {journal} {\bibinfo  {journal}
  {Physical Review A}\ }\textbf {\bibinfo {volume} {46}},\ \bibinfo {pages}
  {6315} (\bibinfo {year} {1992})}\BibitemShut {NoStop}%
\bibitem [{\citenamefont {Meixner}\ \emph {et~al.}(1997)\citenamefont
  {Meixner}, \citenamefont {De~Wit}, \citenamefont {Bose},\ and\ \citenamefont
  {Sch{\"o}ll}}]{Meixner1997J}%
  \BibitemOpen
  \bibfield  {author} {\bibinfo {author} {\bibfnamefont {M.}~\bibnamefont
  {Meixner}}, \bibinfo {author} {\bibfnamefont {A.}~\bibnamefont {De~Wit}},
  \bibinfo {author} {\bibfnamefont {S.}~\bibnamefont {Bose}}, \ and\ \bibinfo
  {author} {\bibfnamefont {E.}~\bibnamefont {Sch{\"o}ll}},\ }\href@noop {}
  {\bibfield  {journal} {\bibinfo  {journal} {Physical Review E}\ }\textbf
  {\bibinfo {volume} {55}},\ \bibinfo {pages} {6690} (\bibinfo {year}
  {1997})}\BibitemShut {NoStop}%
\bibitem [{\citenamefont {Bose}\ \emph {et~al.}(2000)\citenamefont {Bose},
  \citenamefont {Rodin},\ and\ \citenamefont {Sch{\"o}ll}}]{Bose2000J}%
  \BibitemOpen
  \bibfield  {author} {\bibinfo {author} {\bibfnamefont {S.}~\bibnamefont
  {Bose}}, \bibinfo {author} {\bibfnamefont {P.}~\bibnamefont {Rodin}}, \ and\
  \bibinfo {author} {\bibfnamefont {E.}~\bibnamefont {Sch{\"o}ll}},\
  }\href@noop {} {\bibfield  {journal} {\bibinfo  {journal} {Physical Review
  E}\ }\textbf {\bibinfo {volume} {62}},\ \bibinfo {pages} {1778} (\bibinfo
  {year} {2000})}\BibitemShut {NoStop}%
\bibitem [{\citenamefont {Baurmann}\ \emph {et~al.}(2007)\citenamefont
  {Baurmann}, \citenamefont {Gross},\ and\ \citenamefont
  {Feudel}}]{Baurmann2007J}%
  \BibitemOpen
  \bibfield  {author} {\bibinfo {author} {\bibfnamefont {M.}~\bibnamefont
  {Baurmann}}, \bibinfo {author} {\bibfnamefont {T.}~\bibnamefont {Gross}}, \
  and\ \bibinfo {author} {\bibfnamefont {U.}~\bibnamefont {Feudel}},\
  }\href@noop {} {\bibfield  {journal} {\bibinfo  {journal} {Journal of
  Theoretical Biology}\ }\textbf {\bibinfo {volume} {245}},\ \bibinfo {pages}
  {220} (\bibinfo {year} {2007})}\BibitemShut {NoStop}%
\bibitem [{\citenamefont {Just}\ \emph {et~al.}(2001)\citenamefont {Just},
  \citenamefont {Bose}, \citenamefont {Bose}, \citenamefont {Engel},\ and\
  \citenamefont {Sch{\"o}ll}}]{Just2001J}%
  \BibitemOpen
  \bibfield  {author} {\bibinfo {author} {\bibfnamefont {W.}~\bibnamefont
  {Just}}, \bibinfo {author} {\bibfnamefont {M.}~\bibnamefont {Bose}}, \bibinfo
  {author} {\bibfnamefont {S.}~\bibnamefont {Bose}}, \bibinfo {author}
  {\bibfnamefont {H.}~\bibnamefont {Engel}}, \ and\ \bibinfo {author}
  {\bibfnamefont {E.}~\bibnamefont {Sch{\"o}ll}},\ }\href@noop {} {\bibfield
  {journal} {\bibinfo  {journal} {Physical Review E}\ }\textbf {\bibinfo
  {volume} {64}},\ \bibinfo {pages} {026219} (\bibinfo {year}
  {2001})}\BibitemShut {NoStop}%
\bibitem [{\citenamefont {Venkov}\ \emph {et~al.}(2007)\citenamefont {Venkov},
  \citenamefont {Coombes},\ and\ \citenamefont {Matthews}}]{Venkov2007J}%
  \BibitemOpen
  \bibfield  {author} {\bibinfo {author} {\bibfnamefont {N.~A.}\ \bibnamefont
  {Venkov}}, \bibinfo {author} {\bibfnamefont {S.}~\bibnamefont {Coombes}}, \
  and\ \bibinfo {author} {\bibfnamefont {P.~C.}\ \bibnamefont {Matthews}},\
  }\href@noop {} {\bibfield  {journal} {\bibinfo  {journal} {Physica D:
  Nonlinear Phenomena}\ }\textbf {\bibinfo {volume} {232}},\ \bibinfo {pages}
  {1} (\bibinfo {year} {2007})}\BibitemShut {NoStop}%
\bibitem [{\citenamefont {Ledesma-Dur{\'a}n}\ and\ \citenamefont
  {Arag{\'o}n}(2020)}]{Ledesma2020J}%
  \BibitemOpen
  \bibfield  {author} {\bibinfo {author} {\bibfnamefont {A.}~\bibnamefont
  {Ledesma-Dur{\'a}n}}\ and\ \bibinfo {author} {\bibfnamefont {J.~L.}\
  \bibnamefont {Arag{\'o}n}},\ }\href@noop {} {\bibfield  {journal} {\bibinfo
  {journal} {Communications in Nonlinear Science and Numerical Simulation}\
  }\textbf {\bibinfo {volume} {83}},\ \bibinfo {pages} {105145} (\bibinfo
  {year} {2020})}\BibitemShut {NoStop}%
\bibitem [{\citenamefont {Song}\ \emph {et~al.}(2016)\citenamefont {Song},
  \citenamefont {Zhang},\ and\ \citenamefont {Peng}}]{Song2016J}%
  \BibitemOpen
  \bibfield  {author} {\bibinfo {author} {\bibfnamefont {Y.}~\bibnamefont
  {Song}}, \bibinfo {author} {\bibfnamefont {T.}~\bibnamefont {Zhang}}, \ and\
  \bibinfo {author} {\bibfnamefont {Y.}~\bibnamefont {Peng}},\ }\href@noop {}
  {\bibfield  {journal} {\bibinfo  {journal} {Communications in Nonlinear
  Science and Numerical Simulation}\ }\textbf {\bibinfo {volume} {33}},\
  \bibinfo {pages} {229} (\bibinfo {year} {2016})}\BibitemShut {NoStop}%
\bibitem [{\citenamefont {Jiang}\ \emph {et~al.}(2020)\citenamefont {Jiang},
  \citenamefont {An},\ and\ \citenamefont {Shi}}]{Jiang2020J}%
  \BibitemOpen
  \bibfield  {author} {\bibinfo {author} {\bibfnamefont {W.}~\bibnamefont
  {Jiang}}, \bibinfo {author} {\bibfnamefont {Q.}~\bibnamefont {An}}, \ and\
  \bibinfo {author} {\bibfnamefont {J.}~\bibnamefont {Shi}},\ }\href@noop {}
  {\bibfield  {journal} {\bibinfo  {journal} {Journal of Differential
  Equations}\ }\textbf {\bibinfo {volume} {268}},\ \bibinfo {pages} {6067}
  (\bibinfo {year} {2020})}\BibitemShut {NoStop}%
\bibitem [{\citenamefont {Song}\ \emph {et~al.}(2019)\citenamefont {Song},
  \citenamefont {Jiang},\ and\ \citenamefont {Yuan}}]{Song2019J}%
  \BibitemOpen
  \bibfield  {author} {\bibinfo {author} {\bibfnamefont {Y.}~\bibnamefont
  {Song}}, \bibinfo {author} {\bibfnamefont {H.}~\bibnamefont {Jiang}}, \ and\
  \bibinfo {author} {\bibfnamefont {Y.}~\bibnamefont {Yuan}},\ }\href@noop {}
  {\bibfield  {journal} {\bibinfo  {journal} {Journal of Applied Analysis and
  Computation}\ }\textbf {\bibinfo {volume} {9}},\ \bibinfo {pages} {1132}
  (\bibinfo {year} {2019})}\BibitemShut {NoStop}%
\bibitem [{\citenamefont {Wu}\ and\ \citenamefont {Zhao}(2020)}]{Wu2020J}%
  \BibitemOpen
  \bibfield  {author} {\bibinfo {author} {\bibfnamefont {D.}~\bibnamefont
  {Wu}}\ and\ \bibinfo {author} {\bibfnamefont {H.}~\bibnamefont {Zhao}},\
  }\href@noop {} {\bibfield  {journal} {\bibinfo  {journal} {Journal of
  Nonlinear Science}\ }\textbf {\bibinfo {volume} {30}},\ \bibinfo {pages}
  {1015} (\bibinfo {year} {2020})}\BibitemShut {NoStop}%
\bibitem [{\citenamefont {Lv}(2021)}]{Lv2021J}%
  \BibitemOpen
  \bibfield  {author} {\bibinfo {author} {\bibfnamefont {Y.}~\bibnamefont
  {Lv}},\ }\href@noop {} {\bibfield  {journal} {\bibinfo  {journal} {Nonlinear
  Dynamics}\ }\textbf {\bibinfo {volume} {107}},\ \bibinfo {pages} {1357}
  (\bibinfo {year} {2021})}\BibitemShut {NoStop}%
\bibitem [{\citenamefont {Duan}\ \emph {et~al.}(2022)\citenamefont {Duan},
  \citenamefont {Niu},\ and\ \citenamefont {Wei}}]{Duan2022J}%
  \BibitemOpen
  \bibfield  {author} {\bibinfo {author} {\bibfnamefont {D.}~\bibnamefont
  {Duan}}, \bibinfo {author} {\bibfnamefont {B.}~\bibnamefont {Niu}}, \ and\
  \bibinfo {author} {\bibfnamefont {J.}~\bibnamefont {Wei}},\ }\href@noop {}
  {\bibfield  {journal} {\bibinfo  {journal} {Discrete and Continuous Dynamical
  Systems-Series B}\ }\textbf {\bibinfo {volume} {27}},\ \bibinfo {pages}
  {3683} (\bibinfo {year} {2022})}\BibitemShut {NoStop}%
\bibitem [{\citenamefont {Chen}\ \emph {et~al.}(2021)\citenamefont {Chen},
  \citenamefont {Wu}, \citenamefont {Liu},\ and\ \citenamefont
  {Fu}}]{Chen2021J}%
  \BibitemOpen
  \bibfield  {author} {\bibinfo {author} {\bibfnamefont {M.}~\bibnamefont
  {Chen}}, \bibinfo {author} {\bibfnamefont {R.}~\bibnamefont {Wu}}, \bibinfo
  {author} {\bibfnamefont {H.}~\bibnamefont {Liu}}, \ and\ \bibinfo {author}
  {\bibfnamefont {X.}~\bibnamefont {Fu}},\ }\href@noop {} {\bibfield  {journal}
  {\bibinfo  {journal} {Chaos, Solitons \& Fractals}\ }\textbf {\bibinfo
  {volume} {153}},\ \bibinfo {pages} {111509} (\bibinfo {year}
  {2021})}\BibitemShut {NoStop}%
\bibitem [{\citenamefont {Abid}\ \emph {et~al.}(2015)\citenamefont {Abid},
  \citenamefont {Yafia}, \citenamefont {Aziz-Alaoui}, \citenamefont
  {Bouhafaa},\ and\ \citenamefont {Abichoua}}]{Abid2015J}%
  \BibitemOpen
  \bibfield  {author} {\bibinfo {author} {\bibfnamefont {W.}~\bibnamefont
  {Abid}}, \bibinfo {author} {\bibfnamefont {R.}~\bibnamefont {Yafia}},
  \bibinfo {author} {\bibfnamefont {M.}~\bibnamefont {Aziz-Alaoui}}, \bibinfo
  {author} {\bibfnamefont {H.}~\bibnamefont {Bouhafaa}}, \ and\ \bibinfo
  {author} {\bibfnamefont {A.}~\bibnamefont {Abichoua}},\ }\href@noop {}
  {\bibfield  {journal} {\bibinfo  {journal} {Applied Mathematics and
  Computation}\ }\textbf {\bibinfo {volume} {260}},\ \bibinfo {pages} {292}
  (\bibinfo {year} {2015})}\BibitemShut {NoStop}%
\bibitem [{\citenamefont {Paquin-Lefebvre}\ \emph {et~al.}(2019)\citenamefont
  {Paquin-Lefebvre}, \citenamefont {Nagata},\ and\ \citenamefont
  {Ward}}]{Paquin-Lefebvre2019J}%
  \BibitemOpen
  \bibfield  {author} {\bibinfo {author} {\bibfnamefont {F.}~\bibnamefont
  {Paquin-Lefebvre}}, \bibinfo {author} {\bibfnamefont {W.}~\bibnamefont
  {Nagata}}, \ and\ \bibinfo {author} {\bibfnamefont {M.~J.}\ \bibnamefont
  {Ward}},\ }\href@noop {} {\bibfield  {journal} {\bibinfo  {journal} {SIAM
  Journal on Applied Dynamical Systems}\ }\textbf {\bibinfo {volume} {18}},\
  \bibinfo {pages} {1334} (\bibinfo {year} {2019})}\BibitemShut {NoStop}%
\bibitem [{\citenamefont {Golubitsky}\ \emph {et~al.}(1989)\citenamefont
  {Golubitsky}, \citenamefont {Stewart},\ and\ \citenamefont
  {Schaeffer}}]{Golubitsky1989M}%
  \BibitemOpen
  \bibfield  {author} {\bibinfo {author} {\bibfnamefont {M.}~\bibnamefont
  {Golubitsky}}, \bibinfo {author} {\bibfnamefont {I.}~\bibnamefont {Stewart}},
  \ and\ \bibinfo {author} {\bibfnamefont {D.~G.}\ \bibnamefont {Schaeffer}},\
  }\href@noop {} {\emph {\bibinfo {title} {Singularities and {G}roups in
  {B}ifurcation {T}heory: Volume II}}}\ (\bibinfo  {publisher}
  {Springer-Verlag},\ \bibinfo {address} {New {Y}ork},\ \bibinfo {year}
  {1989})\BibitemShut {NoStop}%
\bibitem [{\citenamefont {Guo}\ and\ \citenamefont {Wu}(2013)}]{Guo2013M}%
  \BibitemOpen
  \bibfield  {author} {\bibinfo {author} {\bibfnamefont {S.}~\bibnamefont
  {Guo}}\ and\ \bibinfo {author} {\bibfnamefont {J.}~\bibnamefont {Wu}},\
  }\href@noop {} {\emph {\bibinfo {title} {Bifurcation {T}heory of {F}unctional
  {D}ifferential {E}quations}}}\ (\bibinfo  {publisher} {Springer-Verlag},\
  \bibinfo {address} {New York},\ \bibinfo {year} {2013})\BibitemShut {NoStop}%
\bibitem [{\citenamefont {Hale}(1977)}]{Hale1977M}%
  \BibitemOpen
  \bibfield  {author} {\bibinfo {author} {\bibfnamefont {J.~K.}\ \bibnamefont
  {Hale}},\ }\href@noop {} {\emph {\bibinfo {title} {Theory of {F}unctional
  {D}ifferential {E}quations}}}\ (\bibinfo  {publisher} {Springer-Verlag},\
  \bibinfo {address} {New York},\ \bibinfo {year} {1977})\BibitemShut {NoStop}%
\bibitem [{\citenamefont {Wu}(1996)}]{Wu1996M}%
  \BibitemOpen
  \bibfield  {author} {\bibinfo {author} {\bibfnamefont {J.}~\bibnamefont
  {Wu}},\ }\href@noop {} {\emph {\bibinfo {title} {Theory and {A}pplications of
  {P}artial {F}unctional {D}ifferential {E}quations}}}\ (\bibinfo  {publisher}
  {Springer-Verlag},\ \bibinfo {address} {New {Y}ork},\ \bibinfo {year}
  {1996})\BibitemShut {NoStop}%
\bibitem [{\citenamefont {Murray}(2001)}]{Murray2001M}%
  \BibitemOpen
  \bibfield  {author} {\bibinfo {author} {\bibfnamefont {J.~D.}\ \bibnamefont
  {Murray}},\ }\href@noop {} {\emph {\bibinfo {title} {Mathematical {B}iology
  II: {S}patial {M}odels and {B}iomedical {A}pplications}}}\ (\bibinfo
  {publisher} {Springer-Verlag},\ \bibinfo {address} {New York},\ \bibinfo
  {year} {2001})\BibitemShut {NoStop}%
\bibitem [{\citenamefont {Pinchover}\ and\ \citenamefont
  {Rubinstein}(2005)}]{Pinchover2005M}%
  \BibitemOpen
  \bibfield  {author} {\bibinfo {author} {\bibfnamefont {Y.}~\bibnamefont
  {Pinchover}}\ and\ \bibinfo {author} {\bibfnamefont {J.}~\bibnamefont
  {Rubinstein}},\ }\href@noop {} {\emph {\bibinfo {title} {An {I}ntroduction to
  {P}artial {D}ifferential {E}quations}}}\ (\bibinfo  {publisher} {Cambridge
  University Press},\ \bibinfo {year} {2005})\BibitemShut {NoStop}%
\bibitem [{\citenamefont {Chen}\ \emph {et~al.}(2023)\citenamefont {Chen},
  \citenamefont {Zeng},\ and\ \citenamefont {Niu}}]{Chen2023J}%
  \BibitemOpen
  \bibfield  {author} {\bibinfo {author} {\bibfnamefont {Y.}~\bibnamefont
  {Chen}}, \bibinfo {author} {\bibfnamefont {X.}~\bibnamefont {Zeng}}, \ and\
  \bibinfo {author} {\bibfnamefont {B.}~\bibnamefont {Niu}},\ }\href@noop {}
  {\enquote {\bibinfo {title} {Equivariant {H}opf bifurcation in a class of
  partial functional differential equations on a circular domain},}\ }
  (\bibinfo {year} {2023}),\ \Eprint {http://arxiv.org/abs/2305.05979}
  {arXiv:2305.05979 [math.DS]} \BibitemShut {NoStop}%
\bibitem [{\citenamefont {Guckenheimer}\ and\ \citenamefont
  {Holmes}(1983)}]{Guckenheimer1983M}%
  \BibitemOpen
  \bibfield  {author} {\bibinfo {author} {\bibfnamefont {J.}~\bibnamefont
  {Guckenheimer}}\ and\ \bibinfo {author} {\bibfnamefont {P.}~\bibnamefont
  {Holmes}},\ }\href@noop {} {\emph {\bibinfo {title} {Nonlinear
  {O}scillations, {D}ynamical {S}ystems, and {B}ifurcations of {V}ector
  {F}ields}}}\ (\bibinfo  {publisher} {Springer-Verlag},\ \bibinfo {address}
  {New York},\ \bibinfo {year} {1983})\BibitemShut {NoStop}%
\bibitem [{\citenamefont {Shen}\ and\ \citenamefont {Wei}(2019)}]{Shen2019J}%
  \BibitemOpen
  \bibfield  {author} {\bibinfo {author} {\bibfnamefont {Z.}~\bibnamefont
  {Shen}}\ and\ \bibinfo {author} {\bibfnamefont {J.}~\bibnamefont {Wei}},\
  }\href@noop {} {\bibfield  {journal} {\bibinfo  {journal} {International
  Journal of Bifurcation and Chaos}\ }\textbf {\bibinfo {volume} {29}},\
  \bibinfo {pages} {1950144} (\bibinfo {year} {2019})}\BibitemShut {NoStop}%
\bibitem [{\citenamefont {Steen}(2003)}]{Steen2003J}%
  \BibitemOpen
  \bibfield  {author} {\bibinfo {author} {\bibfnamefont {H.}~\bibnamefont
  {Steen}},\ }\href@noop {} {\bibfield  {journal} {\bibinfo  {journal}
  {Botanica Marina}\ }\textbf {\bibinfo {volume} {46}},\ \bibinfo {pages} {36}
  (\bibinfo {year} {2003})}\BibitemShut {NoStop}%
\bibitem [{\citenamefont {Manoylov}(2009)}]{Manoylov2009J}%
  \BibitemOpen
  \bibfield  {author} {\bibinfo {author} {\bibfnamefont {K.~M.}\ \bibnamefont
  {Manoylov}},\ }\href@noop {} {\bibfield  {journal} {\bibinfo  {journal}
  {Journal of Freshwater Ecology}\ }\textbf {\bibinfo {volume} {24}},\ \bibinfo
  {pages} {145} (\bibinfo {year} {2009})}\BibitemShut {NoStop}%
\bibitem [{\citenamefont {Faria}(2000)}]{Faria2000J}%
  \BibitemOpen
  \bibfield  {author} {\bibinfo {author} {\bibfnamefont {T.}~\bibnamefont
  {Faria}},\ }\href@noop {} {\bibfield  {journal} {\bibinfo  {journal}
  {Transactions of the American Mathematical Society}\ }\textbf {\bibinfo
  {volume} {352}},\ \bibinfo {pages} {2217} (\bibinfo {year}
  {2000})}\BibitemShut {NoStop}%
\bibitem [{\citenamefont {{van}~Gils}\ and\ \citenamefont
  {Mallet-Paret}(1986)}]{Gils1986J}%
  \BibitemOpen
  \bibfield  {author} {\bibinfo {author} {\bibfnamefont {S.~A.}\ \bibnamefont
  {{van}~Gils}}\ and\ \bibinfo {author} {\bibfnamefont {J.}~\bibnamefont
  {Mallet-Paret}},\ }\href@noop {} {\bibfield  {journal} {\bibinfo  {journal}
  {Proceedings of the Royal Society of Edinburgh Section A: Mathematics}\
  }\textbf {\bibinfo {volume} {104}},\ \bibinfo {pages} {279} (\bibinfo {year}
  {1986})}\BibitemShut {NoStop}%
\end{thebibliography}%

\appendix
\section{The proof of Theorem \ref{ET-EH reduce}}\label{thproof}

In this section, we provide the decomposition of the phase space and the derivation of normal forms, by applying the method in \citep{Faria2000J,Song2016J,Jiang2020J}, which leads to the results in Throrem \ref{ET-EH reduce}.

Let $\Lambda_1=\left\{\gamma:\Gamma_{0m_1}(\gamma)=0,~ \mathrm{Re}{\gamma}=0\right\},~\Lambda_2=\left\{\gamma:\tilde{\Gamma}_{nm_2}(\gamma)=0,~ \mathrm{Re}{\gamma}=0\right\}$. Define a bilinear pairing
\begin{equation}\label{bilinear product}
\begin{aligned}
(\psi, \varphi) &=\int_{0}^{R}\int_{0}^{2\pi}r\left[\overline{\psi(0)}\varphi(0)-\int_{-\tau}^{0}\int_{\xi=0}^{\vartheta} \overline{\psi(\xi-\vartheta)}\mathrm{d}\eta(\hat{\nu},\vartheta)\varphi(\xi)  \mathrm{d} \xi\right]\mathrm{d}r\mathrm{d}\theta,~\psi\in\left(\mathscr{C^*}\right)^n,~ \varphi\in\mathscr{C}^n,
\end{aligned}
\end{equation}
where $\left(\mathscr{C^*}\right)^n$ is the dual space of $\mathscr{C}^n$.
By \citep{Hale1977M,Wu1996M}, one can decompose ${C}^n:=C([-1,0],\mathbb{C}^n)$ by $\Lambda_i$ as
$$
{C}^n={P}_i \oplus {Q}_i,~i=1,2,
$$
where ${P}_i$ is the generalised eigenspace associated with $\Lambda_i$ and ${Q}_i=\{\phi \in \mathscr{C}:~\left(\psi,\phi\right)=0,~ for~all~\psi \in {P}_i^*\}$. Here, ${P}_i^*$ is the dual space of ${P}_i$. Suitably, choose the bases $\Phi_{r\theta}^i$ and $\Psi_{r\theta}^i$ of ${P}_i$ and ${P}_i^*$, respectively, such that $\left( \Psi_{r\theta}^i, \Phi_{r\theta}^i \right)=I_{n_i}$, where $n_i=\mathrm{dim}P_i$.
Analogously, the phase space $\mathscr{C}^n$ can be decomposed as
\begin{equation}\label{X}
\mathscr{C}^n=\mathscr{P} \oplus \mathscr{Q},~\mathscr{P}=\mathrm{Im}{\pi},~\mathscr{Q}=\mathrm{Ker}{\pi},
\end{equation}
where $\mathrm{dim}{\mathscr{P}}=\sum_{i=1}^2{n_i}$, and $\pi:~\mathscr{X} \rightarrow \mathscr{P}$ is a projection defined by
\begin{equation}\label{pi}
\pi(U_t)=\sum_{i=1}^2\left(\Phi_{r\theta}^i\langle \Psi_{r\theta}^i,U_t\rangle\right)^\mathrm{T}.
\end{equation}

In Table \ref{assumptionsT}, we list roots with zero real part of the characteristic equation. By (ET-EH), we get that $\Lambda_1=\emptyset,~\Lambda_2=\{\pm\mathrm{i}\omega_{H_3},0\}$.
Let
$$
\begin{aligned}
\Phi_{r\theta}^2=&\left(\Phi^{2(1)}\cdot\hat{\phi}_{nm_2}^c,~\Phi^{2(2)}\cdot\hat{\phi}_{nm_2}^c,~\Phi^{2(3)}\cdot\hat{\phi}_{nm_2}^s,~\Phi^{2(4)}\cdot\hat{\phi}_{nm_2}^s,~\Phi^{2(5)}\cdot\hat{\phi}_{nm_2}^c,~\Phi^{2(6)}\cdot\hat{\phi}_{nm_2}^s\right)\\
=&\left(\xi_H\mathrm{e}^{\mathrm{i}\omega_{H_3}\vartheta}\hat{\phi}_{n_{H_3}m_{H_3}}^c,~\overline{\xi_H}\mathrm{e}^{-\mathrm{i}\omega_{H_3}\vartheta}\hat{\phi}_{n_{H_3}m_{H_3}}^c,~\xi_H\mathrm{e}^{\mathrm{i}\omega_{H_3}\vartheta}\hat{\phi}_{n_{H_3}m_{H_3}}^s,~\overline{\xi_H}\mathrm{e}^{-\mathrm{i}\omega_{H_3}\vartheta}\hat{\phi}_{n_{H_3}m_{H_3}}^s,\right.\\
&\left.~~\xi_T\hat{\phi}_{n_{T_3}m_{T_3}}^c,~{\xi_T}\hat{\phi}_{n_{T_3}m_{T_3}}^s\right),\\
\Psi_{r\theta}^2=&\mathrm{col}\left(\Psi^{2(1)}\cdot\hat{\phi}_{nm_2}^c,~\Psi^{2(2)}\cdot\hat{\phi}_{nm_2}^c,~\Psi^{2(3)}\cdot\hat{\phi}_{nm_2}^s,~\Psi^{2(4)}\cdot\hat{\phi}_{nm_2}^s,~\Psi^{2(5)}\cdot\hat{\phi}_{nm_2}^c,~\Psi^{2(6)}\cdot\hat{\phi}_{nm_2}^s\right)\\
=&\mathrm{col}\left(\eta_H^\mathrm{T}\mathrm{e}^{\mathrm{i}\omega_{H_3}\varrho}\hat{\phi}_{n_{H_3}m_{H_3}}^c,~\overline{\eta_H}^\mathrm{T}\mathrm{e}^{-\mathrm{i}\omega_{H_3}\varrho}\hat{\phi}_{n_{H_3}m_{H_3}}^c,~\eta_H^\mathrm{T}\mathrm{e}^{\mathrm{i}\omega_{H_3}\varrho}\hat{\phi}_{n_{H_3}m_{H_3}}^s,~\overline{\eta_H}^\mathrm{T}\mathrm{e}^{-\mathrm{i}\omega_{H_3}\varrho}\hat{\phi}_{n_{H_3}m_{H_3}}^s\right.\\
&\left.~~~~~~~\eta_T^\mathrm{T}\hat{\phi}_{n_{T_3}m_{T_3}}^c,~\eta_T^\mathrm{T}\hat{\phi}_{n_{T_3}m_{T_3}}^s\right),
\end{aligned}
$$
where $\vartheta \in [-1,0),~\varrho \in (0,1]$, $\xi_H=(p_{11},p_{12},\cdots,p_{1n})^{\mathrm{T}} \in \mathbb{C}^n$ is the eigenvector associated with the eigenvalue $\mathrm{i}\omega$ and $\xi_T=(q_{11},q_{12},\cdots,q_{1n})^{\mathrm{T}} \in \mathbb{R}^n$ is the eigenvector associated with the eigenvalue 0. $\eta_H \in \mathbb{C}^n$ and $\eta_T \in \mathbb{R}^n$ are the corresponding adjoint eigenvectors that satisfy
$
\left( \Psi_{r\theta}^2, \Phi_{r\theta}^2\right)=I_6.
$

According to (\ref{pi}), $U_t=({u_1}_t,{u_2}_t,\cdots,{u_n}_t)$ can be decomposed as
\begin{equation}\label{decompose U}
\begin{aligned}
U_t&=\left(\Phi_{r\theta}^1\langle \Psi_{r\theta}^1,U_t\rangle\right)^\mathrm{T}+\left(\Phi_{r\theta}^2\langle \Psi_{r\theta}^2,U_t\rangle\right)^\mathrm{T}+w_t\\
&=\Phi_{r\theta}^2(z_1,z_2,z_3,z_4,z_5,z_6)^\mathrm{T}+w_t\\
&=\Phi_{r\theta}z+w_t,
\end{aligned}
\end{equation}
with $\Phi_{r\theta}=\Phi_{r\theta}^2,~z=(z_1,z_2,z_3,z_4,z_5,z_6)^\mathrm{T}$, and $w_t \in \mathscr{Q}$. Notice that the part $\Phi_{r\theta}^2(z_1,z_2,z_3,z_4,z_5,z_6)^\mathrm{T}$ stands for the solution on the center manifold, by which solutions on the center manifold are approximatively given.

It is easy to verify that
\begin{equation}\label{Mj1z}
\begin{aligned}
&M_j^1(\mu^l z^p e_k)=D_z(z^p \mu^l e_k)Bz-Bz^p \mu^l e_k=\mathrm{i} \omega_{H_3} \left(p_1-p_2+(-1)^k\right)z^p \mu^l e_k,~ k=1,2,3,4\\
&M_j^1(\mu^l z^p e_k)=D_z(z^p \mu^l e_k)Bz-Bz^p \mu^l e_k=\mathrm{i} \omega_{H_3} \left(p_1-p_2\right)z^p \mu^l e_k,~ k=5,6,
\end{aligned}
\end{equation}
with $M_j^1$ defined in \citep{Faria2000J}, ${B}=\mathrm{diag}\{\mathrm{i}\omega_{H_3},-\mathrm{i}\omega_{H_3},\mathrm{i}\omega_{H_3},-\mathrm{i}\omega_{H_3},0,0\}$, $j \ge 2,~z^p=z_1^{p_1}z_2^{p_2}z_3^{p_3}z_4^{p_4}z_5^{p_5}z_6^{p_6},~\mu^l=\mu_1^{l_1}\mu_2^{l_2},~p_1+p_2+p_3+p_4+l_1+l_2=j$, and $\{e_1,e_2,e_3,e_4,e_5,e_6\}$ being the canonical basis for $\mathscr{C}^6$.
Therefore,
$$
\begin{aligned}
\mathrm{Im}(M_2^1)^{\mathrm{c}}=\mathrm{span}&\left\{z_1z_5e_1,z_2z_5e_2,z_3z_5e_3,z_4z_5e_4,z_3z_5e_1,z_4z_5e_2,z_1z_5e_3,z_2z_5e_4,z_5^2e_5,z_6^2e_5,z_5z_6e_5,\right.\\
&\left.z_1z_6e_1,z_2z_6e_2,z_3z_6e_3,z_4z_6e_4,z_3z_6e_1,z_4z_6e_2,z_1z_6e_3,z_2z_6e_4,z_5^2e_6,z_6^2e_6,z_5z_6e_6,\right.\\
&\left.z_1\mu_ie_1,z_2\mu_ie_2,z_3\mu_ie_3,z_4\mu_ie_4,z_3\mu_ie_1,z_4\mu_ie_2,z_1\mu_ie_3,z_2\mu_ie_4,\right.\\
&\left.z_5\mu_ie_5,z_5\mu_ie_6,z_6\mu_ie_5,z_6\mu_ie_6,\mu_1\mu_2e_5,\mu_i^2e_5,\mu_1\mu_2e_6,\mu_i^2e_6
\right\},\\
\end{aligned}
$$
$$
\begin{aligned}
\mathrm{Im}(M_3^1)^{\mathrm{c}}=\mathrm{span}&\left\{z_1^2z_2e_1,z_3^2z_2e_1,z_1^2z_4e_1,z_3^2z_4e_1,z_1z_5^2e_1,z_1z_6^2e_1,z_1z_5z_6e_1,z_3z_5^2e_1,z_3z_6^2e_1,z_3z_5z_6e_1,\right.\\
&\left.z_1\mu_i^2e_1,z_3\mu_i^2e_1,z_1z_5\mu_ie_1,z_1z_6\mu_ie_1,z_3z_5\mu_ie_1,z_3z_6\mu_ie_1,z_1\mu_1\mu_2e_1,z_3\mu_1\mu_2e_1,\right.\\
&\left.z_1z_2z_3e_1,z_1z_3z_4e_1,\right.\\
&\left.z_1^2z_2e_3,z_3^2z_2e_3,z_1^2z_4e_3,z_3^2z_4e_3,z_1z_5^2e_3,z_1z_6^2e_3,z_1z_5z_6e_3,z_3z_5^2e_3,z_3z_6^2e_3,z_3z_5z_6e_3,\right.\\
&\left.z_1\mu_i^2e_3,z_3\mu_i^2e_3,z_1z_5\mu_ie_3,z_1z_6\mu_ie_3,z_3z_5\mu_ie_3,z_3z_6\mu_ie_3,z_1\mu_1\mu_2e_3,z_3\mu_1\mu_2e_3,\right.\\
&\left.z_1z_2z_3e_3,z_1z_3z_4e_3,\right.\\
&\left.z_2^2z_1e_2,z_4^2z_1e_2,z_2^2z_3e_2,z_4^2z_3e_2,z_2z_5^2e_2,z_2z_6^2e_2,z_2z_5z_6e_2,z_4z_5^2e_2,z_4z_6^2e_2,z_4z_5z_6e_2,\right.\\
&\left.z_2\mu_i^2e_2,z_4\mu_i^2e_2,z_2z_5\mu_ie_2,z_2z_6\mu_ie_2,z_4z_5\mu_ie_2,z_4z_6\mu_ie_2,z_2\mu_1\mu_2e_2,z_4\mu_1\mu_2e_2,\right.\\
&\left.z_1z_2z_4e_2,z_2z_3z_4e_2,\right.\\
&\left.z_2^2z_1e_4,z_4^2z_1e_4,z_2^2z_3e_4,z_4^2z_3e_4,z_2z_5^2e_4,z_2z_6^2e_4,z_2z_5z_6e_4,z_4z_5^2e_4,z_4z_6^2e_4,z_4z_5z_6e_4,\right.\\
&\left.z_2\mu_i^2e_4,z_4\mu_i^2e_4,z_2z_5\mu_ie_4,z_2z_6\mu_ie_4,z_4z_5\mu_ie_4,z_4z_6\mu_ie_4,z_2\mu_1\mu_2e_4,z_4\mu_1\mu_2e_4,\right.\\
&\left.z_1z_2z_4e_4,z_2z_3z_4e_4,\right.\\
&\left.z_1z_2z_5e_5,z_1z_4z_5e_5,z_3z_2z_5e_5,z_3z_4z_5e_5,z_1z_2z_6e_5,z_1z_4z_6e_5,z_3z_2z_6e_5,z_3z_4z_6e_5,\right.\\
&\left.z_1z_2\mu_ie_5,z_1z_4\mu_ie_5,z_3z_2\mu_ie_5,z_3z_4\mu_ie_5,z_5\mu_i^2e_5,z_5\mu_1\mu_2e_5,z_6\mu_i^2e_5,z_6\mu_1\mu_2e_5,\right.\\
&\left.z_5^3e_5,z_6^3e_5,z_5^2z_6e_5,z_6^2z_5e_5,\mu_i^3e_5,\mu_1^2\mu_2e_5,\mu_1\mu_2^2e_5,\right.\\
&\left.z_1z_2z_5e_6,z_1z_4z_5e_6,z_3z_2z_5e_6,z_3z_4z_5e_6,z_1z_2z_6e_6,z_1z_4z_6e_6,z_3z_2z_6e_6,z_3z_4z_6e_6,\right.\\
&\left.z_1z_2\mu_ie_6,z_1z_4\mu_ie_6,z_3z_2\mu_ie_6,z_3z_4\mu_ie_6,z_5\mu_i^2e_6,z_5\mu_1\mu_2e_6,z_6\mu_i^2e_6,z_6\mu_1\mu_2e_6,\right.\\
&\left.z_5^3e_6,z_6^3e_6,z_5^2z_6e_6,z_6^2z_5e_6,\mu_i^3e_6,\mu_1^2\mu_2e_6,\mu_1\mu_2^2e_6
\right\},
\end{aligned}
$$

In fact, according to \citep{Song2016J}, the normal forms for Turing-Hopf bifurcation has the following form
\begin{equation}\label{normal forms for Turing-Hopf}
\dot{z}={B}_1z+\frac{1}{2!}g_2^1(z,0,\mu)+\frac{1}{3!}g_3^1(z,0,0)+o(|z||\mu|^2).
\end{equation}
where $g_2^1(z,0,\mu)$ and $g_3^1(z,0,0)$ are defined in \citep{Song2016J}. By the analysis in \citep{Chen2023J,Faria2000J,Song2016J}, noticing the fact
$$
\int_0^R\int_0^{2\pi} r \left(\hat{\phi}_{n_{H_3}m_{H_3}}^c\right)^{k_1} \left(\hat{\phi}_{n_{H_3}m_{H_3}}^s \right)^{k_2} \left(\hat{\phi}_{n_{T_3}m_{T_3}}^c\right)^{k_3}\left(\hat{\phi}_{n_{T_3}m_{T_3}}^s\right)^{k_4}\mathrm{d}\theta \mathrm{d} r\left\{\begin{array}{cc}
\ne 0, & k_1n_{H_3}-k_2n_{H_3}+k_3n_{T_3}-k_4n_{T_3}=0,\\
=0, & others,
\end{array}\right.
$$
and the relationship of $\Phi_{r\theta}$ and $\Psi_{r\theta}$,
we obtain that when $n_{T_3}\ne 2n_{H_3}$, the normal forms truncated to the third order for ET-EH bifurcation can be summarized as
\begin{equation}\label{ET-EH Normal form}
\begin{aligned}
\dot{z}_1=&\mathrm{i}\omega_{H_3}z_1+B_{11}\mu_1z_1+B_{21}\mu_2z_1+B_{100020}z_1z_5^2+B_{001020}z_3z_5^2+B_{100002}z_1z_6^2+B_{001002}z_3z_6^2\\
&+B_{210000}z_1^2z_2+B_{200100}z_1^2z_4+B_{012000}z_3^2z_2+B_{002100}z_3^2z_4+B_{111000}z_1z_2z_3+B_{101100}z_1z_3z_4\\
&+B_{100011}z_1z_5z_6+B_{001011}z_3z_5z_6,\\
\dot{z}_2=&-\mathrm{i}\omega_{H_3}z_2+\overline{B_{11}}\mu_1z_2+\overline{B_{21}}\mu_2z_2+\overline{B_{100020}}z_2z_6^2+\overline{B_{001020}}z_4z_6^2+\overline{B_{100002}}z_2z_5^2+\overline{B_{001002}}z_4z_5^2\\
&+\overline{B_{210000}}z_1z_2^2+\overline{B_{200100}}z_2^2z_3+\overline{B_{012000}}z_4^2z_1+\overline{B_{002100}}z_4^2z_3+\overline{B_{111000}}z_1z_2z_4+\overline{B_{101100}}z_2z_3z_4,\\
&+\overline{B_{100011}}z_2z_5z_6+\overline{B_{001011}}z_2z_5z_6,\\
\dot{z}_3=&\mathrm{i}\omega_{H_3}z_3+B_{11}\mu_1z_3+B_{21}\mu_2z_3+B_{100020}z_3z_5^2+B_{001020}z_1z_5^2+B_{100002}z_3z_6^2+B_{001002}z_1z_6^2\\
&+B_{210000}z_3^2z_4+B_{200100}z_3^2z_2+B_{012000}z_1^2z_4+B_{002100}z_1^2z_2+B_{111000}z_1z_3z_4+B_{101100}z_1z_2z_3,\\
&+B_{100011}z_3z_5z_6+B_{001011}z_1z_5z_6,\\
\dot{z}_4=&-\mathrm{i}\omega_{H_3}z_4+\overline{B_{11}}\mu_1z_4+\overline{B_{21}}\mu_2z_2+\overline{B_{10002}}z_4z_6^2+\overline{B_{00102}}z_2z_6^2+\overline{B_{100002}}z_4z_5^2+\overline{B_{001002}}z_2z_5^2\\
&+\overline{B_{21000}}z_3z_4^2+\overline{B_{20010}}z_4^2z_1+\overline{B_{01200}}z_2^2z_3+\overline{B_{00210}}z_1^2z_2+\overline{B_{11100}}z_2z_3z_4+\overline{B_{10110}}z_1z_2z_4,\\
&+\overline{B_{100002}}z_4z_6^2+\overline{B_{001002}}z_2z_6^2\\
\dot{z}_5=&B_{15}\mu_1z_5+B_{25}\mu_2z_5+B_{110010}z_1z_2z_5+B_{100110}z_1z_4z_5+B_{011010}z_2z_3z_5+B_{001110}z_3z_4z_5\\
&+B_{110001}z_1z_2z_6+B_{100101}z_1z_4z_6+B_{011001}z_2z_3z_6+B_{001101}z_3z_4z_6\\
&+B_{000030}z_5^3+B_{000003}z_6^3+B_{000021}z_5^2z_6+B_{000012}z_5z_6^2,\\
\dot{z}_6=&B_{15}\mu_1z_6+B_{25}\mu_2z_6+B_{110010}z_1z_2z_6+B_{100110}z_1z_4z_6+B_{011010}z_2z_3z_6+B_{001110}z_3z_4z_6\\
&+B_{110001}z_1z_2z_5+B_{100101}z_1z_4z_5+B_{011001}z_2z_3z_5+B_{001101}z_3z_4z_5\\
&+B_{000030}z_6^3+B_{000003}z_5^3+B_{000021}z_6^2z_5+B_{000012}z_6z_5^2.
\end{aligned}
\end{equation}
By \citep{Gils1986J}, after a sequence of local invertible transformations, the normal form truncated to the third order can be reduced to
\begin{equation}\label{z1z2z3z4z5z6}
\begin{aligned}
&\dot{z}_1=\mathrm{i}\omega_{H_3}z_1+B_{11}\mu_1z_1+B_{21}\mu_2z_1+B_{200100}z_1^2z_4+B_{111000}z_1z_2z_3+B_{100011}z_1z_5z_6,\\
&\dot{z}_2=-\mathrm{i}\omega_{H_3}z_2+\overline{B_{11}}\mu_1z_2+\overline{B_{21}}\mu_2z_2+\overline{B_{200100}}z_3z_2^2+\overline{B_{111000}}z_1z_2z_4+\overline{B_{100011}}z_2z_5z_6,\\
&\dot{z}_3=\mathrm{i}\omega_{H_3}z_3+B_{11}\mu_1z_3+B_{21}\mu_2z_3+B_{200100}z_3^2z_2+B_{111000}z_1z_3z_4+B_{100011}z_3z_5z_6,\\
&\dot{z}_4=-\mathrm{i}\omega_{H_3}z_4+\overline{B_{11}}\mu_1z_4+\overline{B_{21}}\mu_2z_4+\overline{B_{200100}}z_1z_4^2+\overline{B_{111000}}z_2z_3z_4+\overline{B_{100011}}z_4z_5z_6,\\
&\dot{z}_5=B_{15}\mu_1z_5+B_{25}\mu_2z_5+B_{100110}z_1z_4z_5+B_{011010}z_2z_3z_5+B_{000021}z_5^2z_6,\\
&\dot{z}_6=B_{15}\mu_1z_6+B_{25}\mu_2z_6+B_{100110}z_1z_4z_6+B_{011010}z_2z_3z_6+B_{000021}z_5z_6^2.
\end{aligned}
\end{equation}
The proof is similar to Lemma \uppercase\expandafter{\romannumeral3}.2 of \citep{Chen2023J}.

Through the change of variables
\begin{equation}\label{trans}
\begin{array}{lll}
z_1=\rho_{H^1}\mathrm{e}^{\mathrm{i}\chi_{H^1}},~&z_4=\rho_{H^1}\mathrm{e}^{-\mathrm{i}\chi_{H^1}},\\
z_3=\rho_{H^2}\mathrm{e}^{\mathrm{i}\chi_{H^2}},~&z_2=\rho_{H^2}\mathrm{e}^{-\mathrm{i}\chi_{H^2}},\\
z_5=\rho_{T^1},~&z_6=\rho_{T^2},
\end{array}
\end{equation}
we obtain (\ref{rhoETEH}) with
$$
\begin{aligned}
\epsilon_1(\mu)=\mathrm{Re}\{B_{11}\}\mu_1&+\mathrm{Re}\{B_{21}\}\mu_2,~\epsilon_2(\mu)=B_{15}\mu_1+B_{25}\mu_2,\\
c_{11}= \mathrm{Re}\{B_{200100}\},~&c_{12}=\mathrm{Re}\{B_{111000}\},~c_{13}=\mathrm{Re}\{B_{100011}\},\\
c_{21}= B_{100110},~&c_{22}={B_{011010}},~c_{23}={B_{000021}}.
\end{aligned}
$$

\end{document}